\documentclass{amsart}

\numberwithin{equation}{section}

\theoremstyle{plain}
\newtheorem{theorem}{Theorem}[section]
\newtheorem*{theorem*}{Theorem}
\newtheorem{THEOREM}{Theorem}

\newtheorem{proposition}[theorem]{Proposition}
\newtheorem*{proposition*}{Proposition}
\newtheorem{lemma}[theorem]{Lemma}
\newtheorem*{lemma*}{Lemma}
\newtheorem{corollary}[theorem]{Corollary}
\newtheorem*{corollary*}{Corollary}

\newtheorem*{theorem:dPV}{Theorem (dPV)}
\newtheorem*{theorem:dPVI}{Theorem (dPVI)}
\newtheorem*{proposition:isomonodromy}{Proposition (isomonodromy transformation)}

\theoremstyle{definition}
\newtheorem{definition}[theorem]{Definition}
\newtheorem*{definition*}{Definition}

\newtheorem*{notation*}{Notation}

\theoremstyle{remark}
\newtheorem{remark}[theorem]{Remark}
\newtheorem*{remark*}{Remark}

\newtheorem*{remarks*}{Remarks}

\newtheorem*{example*}{Example}

\newcommand\st{\vert}                       
\newcommand\eqd{:=}                     
\newcommand\s[1]{{{\mathrm S}_#1}}                  
\newcommand\p[1]{{\mathbb{P}^{\mathbf{#1}}}}        
\newcommand\A[1]{{\mathbb{A}^{\mathbf{#1}}}}        
\newcommand\iso{{\widetilde\to}}                
\newcommand\isorat{{\widetilde\dashrightarrow}}     
\newcommand\mtop[1]{\mathop{\mathrm{#1}}\nolimits}
\newcommand{\tr}{\mtop{tr}}                 
\newcommand{\im}{\mtop{im}}                 
\newcommand{\coker}{\mtop{coker}} 
\newcommand{\res}{\mtop{res}}                   
\newcommand{\diag}{\mtop{diag}}             

\newcommand\C{\mathbb{C}}
\newcommand\Z{\mathbb{Z}}
\newcommand\V{\mathbb{V}}
\newcommand\PP{\mathbb{P}}

\newcommand\tP{{\widetilde P}}
\newcommand\tK{{\widetilde K}}
\newcommand\tM{{\widetilde M}}
\newcommand\tU{{\widetilde U}}

\newcommand\oM{{\overline M}}

\newcommand\param{\Theta}       
\newcommand\paramv{\Theta^\sharp}   
\newcommand\paramvi{\Theta^\flat}   
\newcommand\vL{{\mathcal L}}
\newcommand\vA{{\mathcal A}}

\newcommand\vR{{\mathcal R}}
\newcommand\vO{{\mathcal O}}
\newcommand\vS{{\mathcal S}}

\newcommand\mA{{A}}
\newcommand\mB{{B}}
\newcommand\mR{{R}}
\newcommand\mS{{S}}
\newcommand\mT{{T}}
\newcommand\mM{{M}}

\author{D.~Arinkin and A.~Borodin}
\title{Moduli spaces of d-connections and difference Painlev\'e equations}

\begin{document}
\begin{abstract}
We show that difference Painlev\'e equations can be interpreted as
isomorphisms of moduli spaces of d-connections on $\p1$ with given
singularity structure. In particular, we derive a difference
equation that lifts to an isomorphism between
$A_2^{(1)*}$--surfaces in Sakai's classification \cite{S}; it
degenerates to both difference Painlev\'e V and classical
(differential) Painlev\'e VI equations. This difference equation
has been known before under the name of asymmetric discrete
Painlev\'e IV equation.
\end{abstract}
\maketitle

\section{Introduction}

This paper is about difference Painlev\'e equations and their
geometric properties. The term `discrete (difference,
$q$-difference, or elliptic) Painlev\'e equation' is rather vague:
there exist different ways of discretizing the classical (2nd
order differential) Painlev\'e equations, see e.g. \cite{GNR},
\cite{ORGT}, \cite{NY}, \cite{JS}, \cite{S}. We will consider the
equations that fit into Sakai's classification described in
\cite{S}.

Any equation of Sakai's hierarchy, by definition, originates from
a birational automorphism of $\C^2$ that lifts to a regular
isomorphism between the blow-ups of $\p2$ at $9$ points.
This geometric property allows to classify the equations
according to the type of the resulting surface. The hierarchy also
includes the classical Painlev\'e equations for which the surfaces
are viewed as spaces of initial conditions, see
\cite{Ok1}, \cite{Ok2}.

In the last few years  several researchers computed the so-called
`gap probabilities' in various
discrete probabilistic models of random matrix type, see
\cite{Ba}, \cite{B1}, \cite{BB}, \cite{AvM},
\cite{FW1}-\cite{FW3}. Surprisingly, these
quantities were often expressible in terms of certain specific
solutions of equations from Sakai's hierarchy. Later it was
demonstrated that the equations arising in probabilistic models
can be viewed as reductions of isomonodromy transformations of
systems of linear difference equations with rational coefficients
\cite{B2}. Further discussion of monodromy for difference equations
can be found in \cite{K}.

The goal of this paper is two-fold. First, we show how the
geometric approach to isomonodromy transformations implies that
the transformations lift to isomorphisms between suitable
surfaces. This provides a conceptual explanation of the above
coincidence. The surfaces are geometrically interpreted as
suitable moduli spaces of {\it d-connections} (short for
``difference connections'') on the Riemann sphere. Second, we
derive an equation of Sakai's hierarchy that lifts to an
isomorphism between $A_2^{(1)*}$--surfaces in Sakai's
classification \cite{S}. We call this equation the {\it difference
Painlev\'e VI}, or dPVI. 

Let us briefly describe our results.

Consider a matrix linear difference equation
\begin{equation}
y(z+1)=A(z)y(z),\qquad A(z)=A_0z^n+\dots+A_{n-1}z+A_n,\quad
A_i\in\operatorname{Mat}(m,\C).
\label{eq:deq}
\end{equation}
We will always assume that $A_0$ is invertible. According to
\cite{B2}, isomonodromy transformations of this equation consist
of maps of the form
\begin{equation}
A(z)\mapsto A'(z)=R(z+1)A(z)R(z)^{-1}
\label{eq:1}
\end{equation}
for suitable rational matrix-valued functions $R(z)$. For generic
$A(z)$, these transformations are parameterized by integral shifts
of the zeros of $A(z)$ and of certain exponents at $z=\infty$ with
total sum of shifts equal to 0, see \cite{B2}[Theorem 2.1].
We can then express the matrix elements of $A'(z)$ as functions
of the matrix elements of $A(z)$; in special cases, the expressions
give rise to the difference Painlev\'e equations.

However, the isomonodromy transformation is defined only when
$A(z)$ is generic enough. Therefore, the resulting maps are
rational rather than regular, that is, the formulas for matrix
elements of $A'(z)$ have singularities. In order to resolve these
singularities it is convenient to use the geometric approach.

Let $\vL$ be a vector bundle on $\p1$ of rank $m$. Assume that we
are given a \emph{d-connection} on $\vL$ which is, by definition, a
linear operator $\vA(z):\vL_z\to \vL_{z+1}$ that depends on $z$
polynomially. (Here $\vL_z$ is the fiber of $\vL$ over $z$.) If $\vL$ is
the trivial vector bundle, $\vA(z)$ is a matrix difference equation
\eqref{eq:deq}.

There is a natural operation on vector bundles with d-connection called
\emph{modification}: it is induced by a rational isomorphism
$\vR:\vL\dashrightarrow \vL'$ between two vector bundles. A
d-connection $\vA$ on $\vL$ then induces a d-connection $\vA'$ on
$\vL'$ and vice versa. Isomonodromy transformations
described above can be viewed as special cases of such
modifications.

Let us consider the example that leads to the difference Painlev\'e
V equation (dPV). Take $m=(\text{rank of  }\vL)=2$; assume that
$\vA(z)$ has four simple zeros $a_1$,$a_2$,$a_3$, $a_4\in\C$, and
that there exists a trivialization of $\vL$ in a neighborhood of
$z=\infty$ with respect to which the matrix of $\vA(z)$ has the
form
$$
A(z)=\begin{bmatrix} \rho_1&0\\0&\rho_2\end{bmatrix} z^2+ \begin{bmatrix}
\rho_1d_1& 0\\0&\rho_2d_2\end{bmatrix} z+O(1), \qquad z\to\infty.
$$

\begin{proposition:isomonodromy} Under
certain non-degeneracy conditions on the parameters
$(a_1,\dots,a_4, \rho_1,\rho_2, d_1,d_2)$, for any vector
bundle $\vL$ with d-connection $\vA$ as above and any integral
shifts of the parameters $a_1,\dots,a_4,d_1,d_2$,
there exists a unique vector bundle $\vL'$ with d-connection $\vA'$
related to $(\vL,\vA)$ by a modification, and such that it
satisfies the above assumptions with shifted values of parameters.
\end{proposition:isomonodromy}

Note that we do not need to assume that $(\vL,\vA)$ is generic.
This means that the modifications of this proposition give (regular,
not birational) isomorphisms of the moduli spaces of vector
bundles with d-connections with given singularity structure, provided
the parameters are generic.

From now on, let us also assume that
$$
\deg(\vL)=-(a_1+\dots+a_4+d_1+d_2)=-1.
$$
This condition implies that $\vL$ is always isomorphic to
$\vO\oplus\vO(-1)$. (Notice that an isomonodromy transformation fixes $\deg(\vL)$
if and only if the corresponding shifts of the parameters
$a_1,\dots,a_4,d_1,d_2$ add up to zero.) By a choice of basis in $\vL$,
the moduli space of d-connections can be
identified with equivalence classes of $2\times 2$ matrices with
polynomial entries satisfying
$$
\gathered
 A=\begin{bmatrix}
a_{11}&a_{12}\\a_{21}&a_{22}\end{bmatrix},\qquad \deg a_{11}\le 2,\
\deg{a_{22}}\le 2,\  \deg a_{21}\le 1,\  \deg a_{12}\le 3,\\
\det A(z)=\operatorname{const}(z-a_1)(z-a_2)(z-a_3)(z-a_4), \\
a_{11}+a_{22}(1+z^{-1})=(\rho_1+\rho_2)z^2+(d_1\rho_1+d_2\rho_2)z+O(1),
\endgathered
$$
modulo the gauge transformations of the form \eqref{eq:1} with
polynomial
$$
R=\begin{bmatrix} r_{11}&r_{12}\\0&r_{22}\end{bmatrix},\qquad
r_{11}=\operatorname{const},\ {r_{22}}=\operatorname{const},\ \deg
r_{12}\le 1.
$$
It is not hard to see that this moduli space is two-dimensional.
We show that its smallest compactification is a surface of the Sakai
type $D_4^{(1)}$; in particular, it is a blow-up of
$\p2$ at 9 points (we use a different description as a blow-up
of $\p1\times\p1$). The moduli space itself is the
complement of $5$ curves (the support of the
unique effective anti-canonical divisor) inside this surface.

In order to connect this picture to dPV, we introduce
coordinates on the moduli spaces.

\begin{theorem:dPV} Take the zero of the linear polynomial
$a_{21}$ as the first coordinate, denote it by $q$, and take the
value of the matrix element $a_{11}$ at $q$ divided by
$(q-a_3)(q-a_4)$ as the second coordinate, denote it by $p$.
Consider the modification of $(\vL,\vA)$ to $(\vL',\vA')$ that shifts
$$
a_1\mapsto a_1-1, \quad a_2\mapsto a_2-1,\quad d_1\mapsto
d_1+1,\quad d_2\mapsto d_2+1.
$$
Then the coordinates $(p',q')$ on the moduli space of $(\vL',\vA')$ are
related to $(p,q)$ by
\begin{equation*}
\begin{cases}
q'+q=a_3+a_4+\dfrac{\rho_1(d_1+a_3+a_4)}{p-\rho_1}+
\dfrac{\rho_2(d_2+a_3+a_4)}{p-\rho_2}\,,
\\
p'p=\dfrac{(q'-a_1+1)(q'-a_2+1)}{(q'-a_3)(q'-a_4)}\cdot\rho_1\rho_2.
\end{cases}
\end{equation*}
\end{theorem:dPV}
This is exactly the dPV equation of \cite{GORS}, \cite{S}.
\footnote{A different reduction of an isomonodromy transformation
to dPV can be found in \cite{B1}.}

\begin{remark*} The idea of using $(q, p)$ as coordinates on the moduli space
is by no means new. For Painlev\'e equations it has been used, for
example, in \cite{JM} in the continuous situation, and in
\cite{JS} in the discrete situation.
\end{remark*}

Another example that we consider in detail deals with rank 2
vector bundles $\vL$ with d-connection $\vA(z)$ which has 6 simple
zeros $a_1,\dots,a_6\in\C$ and whose behavior near $z=\infty$ in a
suitable trivialization is given by
$$
A(z)=\begin{bmatrix} 1&0\\0&1\end{bmatrix} z^3+\begin{bmatrix}
d_1&0\\0&d_2\end{bmatrix} z^2+O(z).
$$
 Quite similarly to the case of dPV discussed
above, there is an action of $\Z^8$ which is parametrized by
integral shifts of $a_i$'s and $d_j$'s. The group
acts by isomorphisms of moduli spaces. Let us again assume
that $\deg(\vL)=-(a_1+\dots+a_6+d_1+d_2)=-1$. Then the corresponding moduli spaces can be
identified with equivalence classes of $2\times 2$ polynomial
matrices satisfying
$$
\gathered
 A=\begin{bmatrix}
a_{11}&a_{12}\\a_{21}&a_{22}\end{bmatrix},\qquad \deg a_{11}\le 3,\
\deg{a_{22}}\le 3,\  \deg a_{21}\le 1,\  \deg a_{12}\le 3,\\
\det A(z)=\operatorname{const} (z-a_1)(z-a_2)(z-a_3)(z-a_4)(z-a_5)(z-a_6)\\
(a_{11}-z^3)(a_{22}(1+z^{-1})- z^3)-a_{12}a_{21}=d_1d_2
z^4+O(z^3),
\endgathered
$$
modulo same gauge transformations as in the case of dPV.

Once again, we show that such a moduli space is two-dimensional
and its smallest compactification can be identified with $\p2$
blown up at $9$ points. The corresponding surface has type
$A_2^{(1)*}$ in Sakai's notation, and the moduli space is the
complement of $3$ curves (the support of the effective anti-canonical divisor) in this surface.

Similarly to the case of dPV, in order to get explicit equations
we need to introduce coordinates on the moduli spaces.

\begin{theorem:dPVI} Take the zero of the matrix element
$a_{21}$ as the first coordinate, denote it by $q$, and take the
value of the matrix element $a_{11}$ at $q$ divided by
$(q-a_4)(q-a_5)(q-a_6)$ as the second coordinate, denote it by
$p$. Consider the modification of $\vL$ to $\vL'$ that shifts
\begin{equation}
a_1\mapsto a_1-1, \quad a_2\mapsto a_2-1,\quad d_1\mapsto
d_1+1,\quad d_2\mapsto d_2+1. \label{eq:2}
\end{equation}
Then the coordinates $(p',q')$ on the moduli space of $\vL'$ are
related to $(p,q)$ by
\begin{equation*}
\begin{cases}
q'=(p-1)(q+1-a_1-a_2)+pa_3+\sum\limits_{j=1,2}\dfrac{c_jp}{\left(q-\frac{p(1-a_1-a_2-d_j)-a_3}{p-1}\right)}
\\
p'\cdot
p=\dfrac{(q'-a_1+1)(q'-a_2+1)}{(q'-a_4)(q'-a_5)(q'-a_6)}\cdot((p-1)(q'-q)+q'-a_3),
\end{cases}
\end{equation*}
where
\begin{equation*}
c_j=\dfrac{(d_j+a_1+a_2+a_4-1)(d_j+a_1+a_2+a_5-1)(d_j+a_1+a_2+a_6-1)}
{(d_j-d_{3-j})} \,.
\end{equation*}
\end{theorem:dPVI}

We will call the relations above the {\it difference Painlev\'e VI
equation\/}.

\begin{remark*} The difference Painlev\'e VI equation above is
equivalent to the asymmetric dPIV equation of \cite{GRO}, see also
earlier references therein. Indeed, introducing the new variable
$r$ instead of $p$ via $p=({q-a_3})/({q+r})$ we can rewrite our
relations as
\begin{equation*}
\begin{cases}
(q+r)(q'+r) =\dfrac{(r+a_3)(r+a_4)(r+a_5)(r+a_6)}
{(r+1-a_1-a_2-d_1)(r+1-a_1-a_2-d_2)}\,,
\\
\\
(q'+r)(q'+r')= \dfrac{(q' - a_3) (q' - a_4) (q' - a_5) (q' -
a_6)}{(q' - (a_1-1)) (q'-(a_2-1))}\,,
\end{cases}
\end{equation*}
which, up to a change of notation, coincides with (1.3) of
\cite{GRO}. We are very grateful to the referee for pointing this
out.

The reason we prefer seemingly more complicated expressions in the
theorem above is that the coordinates have a clear geometric
meaning. This also simplifies various degenerations to other
Painlev\'e equations.
\end{remark*}

It should be noted that formulas for all isomorphisms of Sakai
surfaces in principle can be written using coordinates of
\cite{S}. The computation, however, can be rather tedious.

There are simple degenerations that turn dPVI into dPV and the
classical PVI equations. In a sense, this can be done
simultaneously. Let us consider, in addition to the flow given
by the shift \eqref{eq:2}, the flow generated by the shift
\begin{equation}
a_3\mapsto a_3-1, \quad a_4\mapsto a_4-1,\quad d_1\mapsto
d_1+1,\quad d_2\mapsto d_2+1.
\label{eq:2b}
\end{equation}
Clearly, the flow given by the shift \eqref{eq:2b} is also
described by dPVI with a slightly different $p$-coordinate. Now
let $a_1$, $a_2$, $d_1$, and $d_2$ go to infinity at speeds
$-\rho_1$, $-\rho_2$, $\rho_1$, and $\rho_2$, respectively. In the
limit, the dPVI equation corresponding to \eqref{eq:2} converges
to a continuous vector field which is equivalent to the classical
PVI.\footnote{The classical PVI was also obtained as a limit of
other discrete Painlev\'e equations in \cite{JS}, \cite{RGO}.}
 At the same time, the flow corresponding
to \eqref{eq:2b} converges to a discrete flow described by dPV. As the result, we get two commuting
flows on the same surface (of the Sakai type $D_4^{(1)}$): a vector field given by dPVI and a
discrete dynamics given by dPV.

This limiting picture can be seen from two points of view.
First, the classical PVI possesses the so-called {\it B\"acklund
transformations} which can be described via dPV, see \cite{FW4}.
Second, there is a natural continuous isomonodromy deformation
that moves the parameters $\rho_1,\rho_2$ in the dPV setting; it
can be reduced to the classical PVI. Finally, the geometric
Mellin transform (a version of the Fourier transform) relates the two approaches.
These interrelations (except for the B\"acklund transformations) are
discussed in detail in the body of the paper.

The paper is organized as follows. In Section 1 we state our main
results. In Section 2 we study general properties of d-connections
and discuss various operations on them. Section 3 is dedicated to dPV and the
corresponding moduli space. In Section 4 we describe the relations
between dPV and PVI. Finally, in Section 5 we deal with dPVI, the
associated moduli space, and degenerations of dPVI to dPV and PVI.

The authors are very grateful to D.~Ben-Zvi, P.~Deligne,
D.~Kazhdan, and I.~Krichever for helpful discussions. We would
also like to thank the referees for a number of valuable
suggestions. The first author (D.~A.) was partially supported by
NSF grants DMS-0100108 and DMS-0401164 and by DARPA grant
HR0011-04-1-0031. This research was partially conducted during the
period the second author (A.~B.) served as a Clay Mathematics
Institute Research Fellow. He was also partially supported by the
NSF grant DMS-0402047.

\subsection{Notation}
In this paper, the ground field is $\C$, so `variety' means `variety over $\C$',
`$\p1$' means `$\mathbb{P}^1_\C$' and so on.
$z$ stands for the coordinate on the projective line $\p1$. For a vector bundle
$\vL$ on $\p1$, the fiber of $\vL$ over $z\in\p1$ is denoted by $\vL_z$ and the
space of global sections of $\vL$ is denoted by $\Gamma(\p1,\vL)$.
$\vO(k)$ stands for the line bundle
(vector bundle of rank 1) on $\p1$ whose sections are functions on
$\p1$ with a pole of order at most $k$ (or zero of order at least
$-k$, if $k<0$) at $\infty\in\p1$.

$\diag(\alpha_1,\dots,\alpha_m)$ stands for the diagonal $m\times m$ matrix with entries $\alpha_1,\dots,\alpha_m$.

\section{Main Results}

\subsection{d-connections and their moduli}
Let $\vL$ be a vector bundle on $\p1$ of rank $m$.

\begin{definition} A (rational) \emph{d-connection} on $\vL$ is a linear
operator
\begin{equation*}
\vA(z):\vL_z\to\vL_{z+1}
\end{equation*}
that depends on a point
$z\in\p1-\{\infty\}$ in a rational way (in particular, $\vA(z)$ is
defined for all $z\in\C$ outside of a finite set); here $\vL_z$ is
the fiber of $\vL$ over $z\in\p1$. In other words, $\vA$ is a rational
map between the vector bundle $\vL$ and and its pullback $s^*(\vL)$
via the automorphism $s:\p1\to\p1:z\mapsto z+1$. \label{df:dconn}
\end{definition}

\begin{remark}
Essentially, a d-connection is a system of (rational) linear
difference equations $y(z+1)=\vA(z)y(z)$ on a section $y(z)$ of the
vector bundle $\vL$. Notice that any vector bundle $\vL$ is trivial
when restricted to $\A1=\p1-\{\infty\}$. If we pick a
trivialization $\vS(z):\C^m\iso\vL_z$, $z\in\A1$ of this
restriction (`a basis of $\vL$'), $\vA$ can be written in coordinates
as the matrix-valued function $\mA(z)=\vS(z+1)^{-1}\vA(z)\vS(z)$
(\emph{the matrix of the d-connection}). For two trivializations
$\vS_i(z):\C^m\iso\vL_z$ ($i=1,2$), the corresponding
matrices $\mA_i=\vS_i(z+1)^{-1}\vA(z)\vS_i(z)$ differ by a \emph{d-gauge
transformation}:
\begin{equation*}
\mA_2(z)=\mR(z+1)^{-1}\mA_1(z)\mR(z),
\end{equation*}
for the `d-gauge matrix' $\mR\eqd\vS_1^{-1}\vS_2$. Thus, classification of d-connections is
equivalent to the classification of their matrices up to the d-gauge transformation.
\end{remark}

We work with d-connections that have simple zeroes on $\A1$ and
whose behavior at infinity is `simple' in the sense of the
following definition:

\begin{definition}
Let $\vL$ be a rank 2 vector bundle on $\p1$ and $\vA(z)$ be a d-connection on $\vL$.
Suppose $\vA(z)$ satisfies the following conditions:
\begin{enumerate}
\item The only zeroes and poles of $\vA(z)$ are as follows: a pole
of order $n$ at infinity, and simple zeroes at $k$ distinct points
$a_1,\dots,a_k\in\A1$. Here we say that $a_i$ is a simple zero of
$\vA(z)$ if, at $a_i$,  $\vA(z)$ is regular and $\det(\vA(z))$ has zero
of order $1$.

\item On the formal neighborhood of $\infty\in\p1$, there exist a
trivialization $\vR(z):\C^2\to\vL_z$ ($\vR(z)$ is essentially a
matrix-valued Taylor series in $z^{-1}$) such that the matrix of
$\vA$ with respect to $\vR$ satisfies
\begin{equation}
\vR(z+1)^{-1}\vA(z)\vR(z)=\begin{bmatrix}\rho_1(z^n+d_1z^{n-1})&0\\0&\rho_2(z^n+d_2z^{n-1})\end{bmatrix}
\end{equation}
for some numbers $\rho_1,\rho_2,d_1,d_2\in\C$.
\end{enumerate}

We call such a d-connection $\vA(z)$ (or, more precisely, the pair
$(\vL,\vA)$) a d-connection \emph{of type}
$\theta=(a_1,\dots,a_k;\rho_1,\rho_2,d_1,d_2;n)$.
\label{df:type}
\end{definition}

\begin{remark}
One can also consider d-connections that have simple poles besides
simple zeroes. As it turns out, addition of poles does not lead to
a significantly different object: in Section \ref{sc:operations}, we discuss an
operation (`multiplication by a scalar') that turns a pole of a
d-connection into a zero and vice versa.
\end{remark}

\begin{remark}
The second condition of Definition \ref{df:type} might seem
unnatural, however, Corollary \ref{co:infty} shows that a generic
d-connection satisfies it. See also Remark \ref{rm:formal
solution} for a reformulation of this condition in terms of formal
solutions to a difference equation.
\end{remark}

Denote by $M_\theta$ the moduli space of d-connections of type
$\theta$. One can think of $M_\theta$ in several different ways:
as a set (the set of isomorphism classes of connections of given
type), a category (the category of such connections), a scheme
(the corresponding coarse moduli space), or an algebraic stack
(the fine moduli stack). In this paper, we work with the coarse
moduli space, although some results also hold for other
`incarnations' of $M_\theta$. (Note that we need to impose
some conditions on $\theta$ to make sure that the coarse moduli
space of d-connections of type $\theta$ is a scheme.)

It is easy to see (Corollary \ref{co:degree}) that $M_\theta$ is
empty unless
\begin{gather}
k=2n;\label{cn:pole}\\
\deg(\theta)\eqd-d_1-d_2-\sum_{i=1}^ka_i\mbox{ is an
integer.}\label{cn:degree}
\end{gather}
Let us also consider the following non-degeneracy assumptions on
$\theta$:
\begin{gather}
-d_j-\sum\limits_{i\in I}a_i\not\in\Z\mbox{ for any
}I\subset\{1,\dots,k\}, j=1,2;\label{cn:irred}\\
a_i-a_j\not\in\Z\mbox{ for any }i\ne j;\label{cn:mod}\\
\rho_1,\rho_2\ne0;\rho_1\ne\rho_2.\label{cn:rho}
\end{gather}

Let $\param_{2n}$ be the set of all collections
$\theta=(a_1,\dots,a_{2n};\rho_1,\rho_2,d_1,d_2;n)$, and let
$\paramv_{2n}\subset\param_{2n}$ be the set of $\theta$ that
satisfy \eqref{cn:pole}--\eqref{cn:rho}. Set
$\param=\bigsqcup_n\param_{2n}$,
$\paramv=\bigsqcup_n\paramv_{2n}$.

\begin{remark} Informally speaking, we impose the conditions
\eqref{cn:irred}--\eqref{cn:rho} for the following reasons:
\eqref{cn:mod}, \eqref{cn:rho} simplify modifications of
d-connections (Section \ref{sc:operations}), while \eqref{cn:irred} implies that
d-connections of type $\theta$ are irreducible (Lemma \ref{lm:irred}).
Irreducibility can be used
to prove that the moduli space $M_\theta$ is `nice'; for example,
one can show (using the same ideas as in \cite{AL}) that
$M_\theta$ is a smooth variety of dimension $2n-2$ for any
$\theta\in\paramv_{2n}$. \label{rm:smooth}
\end{remark}

\subsection{Difference $PV$}
We want to study the moduli space $M_\theta$ for
$\theta\in\paramv_{2n}$. As Remark \ref{rm:smooth} shows, the first
interesting case is when $2n=4$: then $M_\theta$ is a smooth
algebraic surface. We will also assume that $\deg(\theta)=-1$ (the
degree is defined in \eqref{cn:degree}).

\begin{remark} The assumption is not too restrictive: using `modifications'
of d-connections (described in Section
\ref{sc:operations}), we can construct for any $\theta$ an isomorphism
$M_\theta\iso M_{\theta'}$, where $\deg(\theta')=-1$.
\end{remark}

We describe the surface $M_\theta$ by introducing
`coordinates' $(q,p)\in(\p1)^2$; more precisely, $M_\theta$ is described
as an open subset in a blow-up of $(\p1)^2$. The construction imitates
the description of the moduli space of connections (\cite{AL}, \cite{IIS2}) which goes
back to Okamoto (\cite{Ok1,Ok2}).

\begin{THEOREM} Suppose
\begin{equation*}
\theta=(a_1,a_2,a_3,a_4;\rho_1,\rho_2,d_1,d_2;2)\in\paramv_4
\end{equation*}
has $\deg(\theta)=-1$. Let $\sigma_1:K_1\to(\p1)^2$ be the blow-up
of $(\p1)^2$ at the following $6$ points: $(q,p)=(a_1,0)$,
$(a_2,0)$, $(a_3,\infty)$, $(a_4,\infty)$, $(\infty,\rho_1)$, and
$(\infty,\rho_2)$ (here $q$ and $p$ are the projections
$(\p1)^2\to\p1$). Consider the two exceptional curves
$E_j=\sigma_1^{-1}(\infty,\rho_j)\subset K_1$, $j=1,2$;
homogeneous coordinates on $E_j$ are given by $(1/q:p-\rho_j)$.
Let $\sigma_2:K_2\to K_1$ be the blow-up of $K_1$ at the two
points $(1/q:p-\rho_j)=(1:\rho_1(d_1+a_3+a_4))$, $j=1,2$ (one
point on each exceptional curve).
\begin{enumerate}
\item There exists an open embedding
$P_2:M_\theta\hookrightarrow K_2$.

\item The complement to $P_2(M_\theta)$ in $K_2$ is the union of
the proper preimages of the following curves:
$\p1\times\{0\},\p1\times\{\infty\},\{\infty\}\times\p1\subset(\p1)^2$,
and the two exceptional curves $E_j\subset K_1$, $j=1,2$.
\end{enumerate}
\label{th:Geometry of dPV}
\end{THEOREM}

\begin{remark} $K_2$ is the
smallest smooth compactification of $M_\theta$ (cf.
\cite[Corollary 5]{AL}): any open embedding
$M_\theta\hookrightarrow\oM$ with smooth projective $\oM$ induces
a regular morphism $\oM\to K_2$. Note also that
$(K_2,K_2-M_\theta)$ is an Okamoto-Painlev\'e pair (of type
$\tilde D_4$) in the sense of \cite{ST}, \cite{STT}; in
particular, $K_2$ is a surface of the Sakai type $D^{(1)}_4$.
\label{rm:smallest}
\end{remark}

In particular, the composition $P:M_\theta\hookrightarrow
K_2\to(\p1)^2$ is birational. Therefore, one can view the
components of $P$ is a kind of `rational coordinates' on
$M_\theta$. We denote the components by $q$ and $p$, so that
$P=(q,p)$.

The natural operations on d-connections (modifications and
multiplications by scalar) define isomorphisms between the spaces
$M_\theta$ for different $\theta$ (see Section \ref{sc:operations}). Our
next result describes such an isomorphism for one of the simplest
modifications of d-connections. The description can be viewed as a
non-linear difference equation in `coordinates' $(p,q)$ (the difference $PV$).

As before, suppose
\begin{equation*}
\theta=(a_1,a_2,a_3,a_4;\rho_1,\rho_2,d_1,d_2;2)\in\paramv_4
\end{equation*}
has $\deg(\theta)=-1$. Set
\begin{equation*}
\theta'=(a_1-1,a_2-1,a_3,a_4;\rho_1,\rho_2,d_1+1,d_2+1;2)\in\paramv_4.
\end{equation*}
Modification of d-connections defines an isomorphism
$dPV:M_\theta\to M_{\theta'}$. Explicitly, for every $(\vL,\vA)\in
M_\theta$, the image $dPV(\vL,\vA)=(\vL',\vA')$ is the only d-connection
of type $\theta'$ that admits a rational isomorphism $\vR:\vL'\isorat
\vL$ that agrees with the d-connections: $\vR(z+1)\vA'(z)=\vA(z)\vR(z)$.

\begin{THEOREM} Set $p'\eqd p\circ dPV,q'\eqd q\circ dPV:M_\theta\to\p1$.
Then
\begin{equation}
\begin{cases}
q'+q=a_3+a_4+\dfrac{\rho_1(d_1+a_3+a_4)}{p-\rho_1}+\dfrac{\rho_2(d_2+a_3+a_4)}{p-\rho_2}
\cr p'\cdot
p=\dfrac{(q'-a_1+1)(q'-a_2+1)}{(q'-a_3)(q'-a_4)}\cdot\rho_1\rho_2
\end{cases}
\label{eq:dPV}
\end{equation}
\label{th:dPV}
\end{THEOREM}

\subsection{Difference $PV$ and classical $PVI$}
As we mentioned above, d-connections and ordinary connections have many common properties. Let us consider
the following class of (ordinary) connections:

Denote by $\Lambda\subset\C^8$ the set of all collections
$\lambda=(\lambda^-_1,\lambda^+_1,\dots,\lambda^-_4,\lambda^+_4)$
such that
\begin{equation*}
\sum_{i=1}^4(\lambda_i^- + \lambda_i^+)\in\Z,\qquad
\lambda_i^+-\lambda_i^-\not\in\Z,\qquad
\sum_{i=1}^4\lambda_i^{\epsilon_i}\not\in\Z
\end{equation*}
for any choice of upper indexes $\epsilon_i\in\{+,-\}$. Let $X\subset(\p1)^4$ be the set of all collections
$x=(x_1,\dots,x_4)$ of four
distinct points of $\p1$:
\begin{equation*}
X\eqd\{(x_1,\dots,x_4)\st x_i\ne x_j\text{ for }i\ne j\}\subset(\p1)^4.
\end{equation*}

\begin{definition} Suppose $(x,\lambda)\in X\times\Lambda$. A
\emph{connection of type $(x,\lambda)$} is a pair $(\vL,\nabla)$ such that $\vL$ is
a rank 2 vector bundle on $\p1$, $\nabla:\vL\to\vL\otimes\Omega_\p1(x_1+\dots+x_4)$ is a connection with simple
poles at $x_i$'s, and the residue of $\nabla$ at $x_i$ has eigenvalues $\{\lambda_i^-,\lambda_i^+\}$.
\end{definition}

For $(x,\lambda)\in X\times\Lambda$, we denote the coarse moduli space of connections of type $(x,\lambda)$ by
$M_{(x,\lambda)}$.
It can be thought of as the space of initial conditions of the Painlev\'e equation $PVI$. The space
$M_{(x,\lambda)}$ has a geometric descritption, which goes back to K.~Okamoto; we remind the description
in Proposition \ref{pp:Geometry of PVI}.
It is easy to see from the description that $M_\theta$ and $M_{(x,\lambda)}$ are isomorphic
for a suitable choice of parameters (they are both surfaces of type $D^{(1)}_4$):

\begin{THEOREM} Suppose
\begin{equation*}
\theta=(a_1,a_2,a_3,a_4;\rho_1,\rho_2,d_1,d_2;2)\in\paramv_4
\end{equation*}
has $\deg(\theta)=-1$. Set
\begin{equation*}
x=(x_1,x_2,x_3,x_4)\eqd(0,\rho_1,\rho_2,\infty)\in X,
\end{equation*}
\begin{equation*}
\lambda=(\lambda_1^-,\lambda_1^+,\dots,\lambda_4^-,\lambda_4^+)\eqd
(a_1,a_2,0,d_1+a_3+a_4,0,d_2+a_3+a_4,-a_3,-a_4)\in\Lambda.
\end{equation*}
Then $M_\theta\simeq M_{(x,\lambda)}$.
\label{th:dPV and PVI}
\end{THEOREM}

\begin{remark}
Theorem \ref{th:dPV and PVI} can be proved by direct calculations, but it can also be explained
in terms of moduli spaces. In Section \ref{sc:Fourier}, we describe
a one-to-one correspondence between d-connections
of type $\theta$ and connections of type $(x,\lambda)$. Up to small `twists', the correspondence
is the geometric Mellin transform of \cite{L}; it is constructed using de Rham cohomology and equivariant cohomology
groups. The Mellin transform is a particular case of the duality for generalized one-motives
(also defined in \cite{L}).
\end{remark}

Now let us fix $\lambda\in\Lambda$ and consider surfaces $M_{(x,\lambda)}$ for all $x\in X$.
They can be viewed as fibers of an algebraic family $M_\lambda\to X$.
The sixth Painlev\'e equation $PVI$ is an algebraic connection on this
family; the (analytic) integral curves of $PVI$ correspond to isomonodromy deformation of connections.

By Theorem \ref{th:dPV and PVI}, the sixth Painlev\'e equation $PVI$
induces a connection on a family of moduli spaces of d-connections.
It turns out that this connection can be defined for arbitrary $\theta\in\paramv_{2n}$ (not necessarily when $2n=4$).
More precisely:

\begin{THEOREM} \label{th:isomonodromic}
Let $n$ be a positive integer.
Fix $a_1,\dots,a_{2n},d_1, d_2\in\C$ that satisfy \eqref{cn:degree}--\eqref{cn:mod}, and set
$P=\{(\rho_1,\rho_2)\in\C^2:\rho_1,\rho_2\ne0,\rho_1\ne\rho_2\}$. For all $\rho\eqd(\rho_1,\rho_2)\in P$, set
$\theta(\rho)=(a_1,\dots,a_{2n};\rho_1,\rho_2,d_1,d_2;n)\in\paramv_{2n}$, and consider the coarse
moduli spaces $M_{\theta(\rho)}$. Clearly, they form a family $M\to P$.
\begin{enumerate}
\item The family
$M\to P$ carries a natural algebraic connection (defined in Section \ref{sc:isomonodromic}).
\item
In the case $2n=4$, this connection coincides with the $PVI$ connection under the isomorphism of
Theorem \ref{th:dPV and PVI}.
\end{enumerate}
\end{THEOREM}

\begin{remark} The connection of Theorem \ref{th:isomonodromic} can be thought of as a `continuous' isomonodromy
deformation of d-connections.
\end{remark}

\subsection{Difference $PVI$}
So far, we have worked with d-connections of type $\theta$, where $\theta$ is non-degenerate
in the sense of \eqref{cn:irred}--\eqref{cn:rho}. It turns out that a different
class of d-connections enjoys similar properties. Namely, let us
replace \eqref{cn:rho} with the following condition:

\begin{equation}
\rho_1=\rho_2\ne0;d_1\ne d_2.\label{cn:rhoalt}
\end{equation}

Let $\paramvi_{2n}\subset\param_{2n}$ be the set of all $\theta$ that satisfy \eqref{cn:pole}--\eqref{cn:mod},
and \eqref{cn:rhoalt}, and set $\paramvi=\bigsqcup_n\paramvi_{2n}$. It can be shown that for
$\theta\in\paramvi_{2n}$, the coarse moduli space $M_\theta$ is a smooth variety of dimension $2n-4$ (recall
that for $\theta\in\paramv_{2n}$, we have $\dim(M_\theta)=2n-2$). Therefore, the first `interesting' case is
$\theta\in\paramvi_6$; then $M_\theta$ is an algebraic surface. As before, we assume $\deg(\theta)=-1$.

Similarly to Theorem \ref{th:Geometry of dPV}, we can describe the moduli space $M_\theta$ using
`coordinates' $(q,p)\in(\p1)^2$.

\begin{THEOREM} Suppose
\begin{equation*}
\theta=(a_1,a_2,a_3,a_4,a_5,a_6;\rho,\rho,d_1,d_2;3)\in\paramvi_6
\end{equation*}
has $\deg(\theta)=-1$. Let $\sigma_1:K_1\to(\p1)^2$ be the blow-up
of $(\p1)^2$ at the following $7$ points: $(q,p)=(a_1,0)$,
$(a_2,0)$, $(a_3,0)$, $(a_4,\infty)$, $(a_5,\infty)$, $(a_6,\infty)$, and
$(\infty,\rho)$ (here $q$ and $p$ are the projections
$(\p1)^2\to\p1$). Consider the exceptional curve
$E=\sigma_1^{-1}(\infty,\rho)\subset K_1$;
a homogeneous coordinate on $E$ is given by $(1/q:p-\rho)$.
Let $\sigma_2:K_2\to K_1$ be the blow-up of $K_1$ at the two
points $(1/q:p-\rho)=(1:\rho(d_j+a_4+a_5+a_6))$, $j=1,2$.
\begin{enumerate}
\item There exists an open embedding
$P_2:M_\theta\hookrightarrow K_2$.

\item The complement to $P_2(M_\theta)$ in $K_2$ is the union of
the proper preimages of the following curves:
$\p1\times\{0\},\p1\times\{\infty\},\{\infty\}\times\p1\subset(\p1)^2$,
and the exceptional curve $E\subset K_1$, $j=1,2$.
\end{enumerate}
\label{th:Geometry of dPVI}
\end{THEOREM}

\begin{remark} Using multiplication by scalar, it is easy to see that the moduli space
$M_\theta$ for $\theta=(a_1,\dots,a_{2n};\rho,\rho,d_1,d_2;n)$ does not depend on $\rho$. Therefore,
we can assume that $\rho=1$ without loss of generality.
\end{remark}

\begin{remark} $K_2$ is not the
smallest smooth compactification of $M_\theta$ (unlike the case when $\theta\in\paramv_4$, see Remark
\ref{rm:smallest}). Indeed, the proper preimage of $\{\infty\}\times\p1\subset(\p1)^2$ is an exceptional curve
in $K_2-M_\theta$. Contracting the exceptional curve, we obtain the smallest smooth compactification of $M_\theta$,
which is a surface of the Sakai type $A_2^{(1)*}$.
\end{remark}

Modifications of d-connections define natural isomorphisms between spaces $M_\theta$.
Similarly to Theorem \ref{th:dPV}, we describe a simple isomorphism of this kind explicitly.
We call the resulting difference equation `the difference $PVI$': as we will see, it degenerates
into both the difference $PV$ (Section \ref{sc:dPVI and dPV}) and the usual $PVI$ (Section \ref{sc:dPVI and PVI}).

Suppose
\begin{equation*}
\theta=(a_1,a_2,a_3,a_4,a_5,a_6;1,1,d_1,d_2;3)\in\paramvi_6
\end{equation*}
has $\deg(\theta)=-1$. Set
\begin{equation*}
\theta'=(a_1-1,a_2-1,a_3,a_4,a_5,a_6;1,1,d_1+1,d_2+1;3)\in\paramvi_6.
\end{equation*}
Modification of d-connections induces an isomorphism
$dPVI:M_\theta\to M_{\theta'}$. Explicitly, for every $(\vL,\vA)\in
M_\theta$, the image $dPVI(\vL,\vA)=(\vL',\vA')$ is the only d-connection
of type $\theta'$ that admits a rational isomorphism $\vR:\vL'\isorat
\vL$ that agrees with the d-connections: $\vR(z+1)\vA'(z)=\vA(z)\vR(z)$.

\begin{THEOREM} Set $p'\eqd p\circ dPVI,q'\eqd q\circ dPVI:M_\theta\to\p1$.
For $j=1,2$, set
\begin{equation*}
c_j\eqd\dfrac{(d_j+a_1+a_2+a_4-1)(d_j+a_1+a_2+a_5-1)(d_j+a_1+a_2+a_6-1)}{(d_j-d_{3-j})}
\end{equation*}
(the denominator is $\pm(d_1-d_2)$).
Then
\begin{equation}
\begin{cases}
q'=(p-1)(q+1-a_1-a_2)+pa_3+\sum_{j=1,2}\dfrac{c_jp}{\left(q-\frac{p(1-a_1-a_2-d_j)-a_3}{p-1}\right)}

\cr
p'\cdot p=\dfrac{(q'-a_1+1)(q'-a_2+1)}{(q'-a_4)(q'-a_5)(q'-a_6)}\cdot((p-1)(q'-q)+q'-a_3)
\end{cases}
\label{eq:dPVI}
\end{equation}
\label{th:dPVI}
\end{THEOREM}

\begin{remark}
Theorem \ref{th:dPV and PVI} identifies $M_\theta$ for $\theta\in\paramv_4$ with a moduli space of connections
of certain kind. A similar statement holds for
$
\theta=(a_1,\dots,a_6;\rho,\rho,d_1,d_2;3)\in\paramvi_6$. In this case,
$M_\theta$ is isomorphic to the moduli space of pairs $(\vL,\nabla)$, where $\vL$ is a rank $3$ bundle on $\p1$
and $\nabla$ is a connection on $\vL$ with first order poles at $\rho$, $0$, and $\infty$ (and no other poles);
the residues at the poles have eigenvalues $\{0,d_1+a_4+a_5+a_6,d_2+a_4+a_5+a_6\}$, $\{a_1,a_2,a_3\}$, $\{-a_4,-a_5,-a_6\}$,
respectively. The isomorphism can be constructed using the Mellin transform (similarly to Section \ref{sc:Fourier}).

Notice that if we interpret $M_\theta$ as a moduli space of rank $3$ bundles with connections on $\p1$, then dPVI
becomes an isomorphism between such moduli spaces (a B\"acklund transformation) which corresponds to
a modification of such bundles.
\end{remark}

\section{General d-connections}
\label{sc:General d-connections}

\subsection{Formal behavior at infinity}
Let $\vL$ be a vector bundle on $\p1$ and $\vA(z):\vL_z\to\vL_{z+1}$ be a
rational d-connection on $\vL$. Since $\infty\in\p1$ is the only
fixed point of the transformation $z\mapsto z+1$, it is natural to
study the restriction of $\vA$ to a neighborhood of infinity. Here
the word `neighborhood' can be understood either analytically (a
small disk) or formally (the formal disk). In this section, we
work with the formal neighborhood: the corresponding
classification problem is significantly easier. The situation is
somewhat similar to classification of irregular singularities for
ordinary differential equations: the formal classification is much
simpler than the analytic one (because of the Stokes' phenomenon).

In the language of difference equations, the problem is to
classify matrices $\mA(z)$ over the ring of formal Laurent series
$\C((z^{-1}))$ modulo d-gauge transformations
\begin{equation*}
\mA(z)\mapsto\mR(z+1)^{-1}\mA(z)\mR(z),
\end{equation*}
where the gauge matrix $\mR(z)$ is an invertible matrix over the
ring of formal Taylor series $\C[[z^{-1}]]$.

If $\mA$ is generic, the answer is given by the following easy
statement (see for instance \cite[Proposition 1.1]{B2}):

\begin{proposition} Suppose that the $m\times m$ matrix
$\mA(z)=\sum\limits_{i\le n}\mA_i z^i$ over $\C((z^{-1}))$
satisfies the following condition:
\begin{equation}
\label{cn:dpv} \begin{array}{l} \mbox{All eigenvalues of the
leading term }\mA_n\mbox{ are distinct and non-zero}\cr\mbox{(in
other words, $\mA_n$ is invertible, regular, and
semisimple).}\end{array}
\end{equation}
Then there exists a gauge matrix
$\mR(z)=\sum\limits_{i\le0}\mR_iz^i$ with invertible $\mR_0$ such
that
\begin{equation}
\mR(z+1)^{-1}\mA(z)\mR(z)=\mA'_n z^n+\mA'_{n-1} z^{n-1},
\end{equation}
where $\mA'_n$ and $\mA'_{n-1}$ are diagonal matrices. The matrix
$\mR(z)$ is uniquely determined up to right multiplication by a
permutation matrix and a constant diagonal matrix. \qed
\label{pp:infty}
\end{proposition}

Denote the diagonal entries of $\mA'_n$ by
$\rho_1,\dots,\rho_m$; notice that $\rho_i$'s are the eigenvalues
of $\mA_n$, in particular, all $\rho_i$ are distinct and non-zero.
Denote the corresponding diagonal entries of $\mA'_{n-1}$ by
$c_1,\dots,c_m$. Set $d_i\eqd c_i/\rho_i$; we
work with $d_i$ rather then $c_i$ because it
simplifies formulas \eqref{cn:degree}, \eqref{cn:irred}. We call
the collection $(\rho_1,\dots,\rho_m,d_1,\dots,d_m;n)$ \emph{the
formal type} of $\mA(z)$ at infinity. Proposition \ref{pp:infty}
implies that the formal type is determined by $\mA(z)$ up to a
simultaneous permutation of $\rho_i$'s and $d_i$'s, that is, up to
the action of the symmetric group $\s{m}$.

\begin{remark}\label{rm:formal solution}
Proposition \ref{pp:infty} is sometimes (for instance, in
\cite{B2}) formulated in terms of formal solutions to the
difference equation: the claim is that the equation
$Y(z+1)=\mA(z)Y(z)$ has a formal solution of the form
\begin{equation*}
Y(z)=(\Gamma(z))^n\left(\sum_{i\le0}{\hat
Y}_iz^i\right)\diag(\rho_1^z z^{d_1},\dots,\rho_m^zz^{d_m}),
\end{equation*}
where $\hat Y_i$ are $m\times m$ matrices, $\hat Y_0$ is
invertible, and $\rho_1,\dots,\rho_m,d_1,\dots,d_m\in\C$.

Note that $\sum_{i\le0}{\hat Y}_iz^i$ does not coincide with
$R(z)$ of Proposition \ref{pp:infty}.
\end{remark}

\begin{remark} The formal type of $\mA(z)$ can be determined
directly, without diagonalizing $\mA(z)$. Indeed, denote by
$\sigma_i(z)$ and $\sigma'_i(z)$ ($i=1,\dots,m$) the coefficients
of the characteristic polynomials of $\mA(z)$ and
$\mR(z+1)^{-1}\mA(z)\mR(z)$ respectively, so that
$\sigma_1(z)=-\tr\mA(z)$ and $\sigma_m(z)=(-1)^m\det\mA(z)$.
Clearly, $\sigma_i(z)$ and $\sigma'_i(z)$ have pole of order
$i\cdot n$ at infinity. One can easily check that the order of
pole of $\sigma_i(z)-\sigma'_i(z)$ is at most $i\cdot n-2$. Thus,
the two leading terms of $\sigma_i(z)$ and $\sigma'_i(z)$
coincide. It is now easy to see that the formal type of $\mA(z)$ is determined
(up to the $\s{m}$-action) by the pairs of leading terms of $\sigma_i(z)$,
$i=1,\dots,m$.

In particular, if we assume $\mA_n$ is diagonal, then its diagonal
entries are the $\rho_i$'s, and the diagonal entries of
$\mA_{n-1}$ equal $\rho_i d_i$, even if $\mA_{n-1}$ is not
diagonal.
\label{rm:characteristic polynomial}
\end{remark}

Let us now translate Proposition \ref{pp:infty} into the language
of d-connections. For simplicity, we only consider vector bundles
of rank $2$.

\begin{corollary}
\label{co:infty} Let $\vA(z)$ be a d-connection on a rank $2$ vector
bundle $\vL$. Denote by $n$ the order of pole of $\vA$ at infinity and
by $\vA_n:\vL_\infty\to\vL_\infty$ the leading term of $\vA$ (that is,
$n$ is the smallest number such that the limit
\begin{equation*}
\vA_n\eqd\lim_{z\to\infty}\vA(z)z^{-n}
\end{equation*}
exists). Suppose all eigenvalues of $\vA_n$ are distinct and
non-zero. Then $\vA(z)$ satisfies the second condition of Definition
\ref{df:type} (for some $\rho_1,\rho_2,d_1,d_2\in\C$). \qed
\end{corollary}

We call the collection $(\rho_1,\rho_2,d_1,d_2;n)$ \emph{the
formal type} of the d-connection $\vA(z)$. It is determined by
$\vA(z)$ up to the action of $\s{2}$. Notice also that in the
situation of Corollary \ref{co:infty}, the condition
\eqref{cn:rho} holds automatically.

\subsection{Operations on d-connections}
\label{sc:operations}
Let us now discuss some natural operations on d-connections. As we
will see, the operations allow us to identify the moduli spaces
(or moduli stacks, or sets of isomorphism classes, or categories) of
d-connections of type $\theta$ for different $\theta$. As a
trivial example, notice that $M_{\theta'}=M_\theta$ if $\theta'$
is obtained from $\theta$ by a permutation of $a_i$'s or a
simultaneous permutation of $\rho_i$'s and $d_i$'s.

{\bf Multiplication by a scalar:} Let $f(z)\ne 0$ be a rational
function on $\p1$, and let $\vA(z)$ be a d-connection on a vector
bundle $\vL$. Clearly, the product $f(z)\vA(z)$ is again a
d-connection on $\vL$.

In the language of difference equations, this operation corresponds to multiplication of
solutions by $\Gamma$-functions. Indeed, let us write
$f(z)=c\prod(z-z_i)^{k_i}$.
Then $y(z)$ solves the difference equation $y(z+1)=\vA(z)y(z)$
if and only if $\tilde y(z)=e^{cz}\prod\Gamma(z-z_i)^{k_i}y(z)$ solves
$\tilde y(z+1)=(f(z)\vA(z))\tilde y(z)$.

On the other hand, multiplication by a scalar is also a special case of a
tensor product of d-connections. We can view $f(z)$ as a
d-connection on the trivial rank one bundle $\vO_\p1$; then
$f(z)\vA(z)$ becomes the natural d-connection on the tensor product
$\vL=\vL\otimes\vO_\p1$ of two vector bundles with d-connections.

\begin{remark} For any d-connection $\vA(z)$, we can pick a function $f(z)$ so that the only pole of the product
$f(z)\vA(z)$ is at infinity. For instance, suppose $\vL$ has rank $2$, and the d-connection
$\vA(z)$ has a simple pole at $z=z_0$; this means that all matrix elements of $A(z)$ (in some basis)
have at most a simple pole and $\det(A(z))$ has a simple pole at $z=z_0$. Then $(z-z_0)A(z)$ has a simple zero at $z_0$.
In this way, classification of rank 2 d-connections with simple poles and simple zeroes on $\p1-\{\infty\}$
is reduced to classification of d-connections with simple zeroes only.
\label{rm:poles and zeroes}
\end{remark}

Now suppose $(\vL,\vA(z))\in M_\theta$ for $\theta\in\param$.
Let $f(z)$ be a rational function; clearly, the product
$(\vL,f(z)\vA(z))$ is a d-connection of type $\theta'$ (for some
$\theta'\in\param$) if and only if the function $f(z)=c$
is a non-zero constant. If $f(z)=c\in\C-\{0\}$, then $(\vL,c\vA)\in
M_{\theta'}$ for
\begin{equation*}
\theta'=(a_1,\dots,a_k;c\rho_1,c\rho_2,d_1,d_2;n).
\end{equation*}
Clearly, the correspondence $(\vL,\vA)\mapsto (\vL,c\vA)$ gives an
isomorphism $\mu=\mu_c:M_\theta\iso M_{\theta'}$; the inverse map
is $\mu_{c^{-1}}$.

{\bf Modification:}
Suppose $\vR:\vL\isorat\vL'$ is a rational isomorphism between two vector bundles $\vL$ and $\vL'$ on $\p1$. Then
a d-connection $\vA(z)$ on $\vL$
induces a d-connection $\vA'$ on $\vL'$ (and vice versa).

In the language of difference equation, this operation is the
d-gauge transformation
\begin{equation}
\mA'(z)=\mR(z+1)^{-1}\mA(z)\mR(z), \label{eq:gauge}
\end{equation}
where $\mR$, $\mA$, and $\mA'$ are the matrices of $\vR$, $\vA$, and
$\vA'$ respectively (corresponding to some choice of bases). We call
$\vA'$ a modification of $\vA$ (of course, $\vA$ is also a
modification of $\vA'$).

\begin{remark} Modifications can also be viewed as an isomonodromy deformation in the sense of \cite{B2}.
Indeed, the monodromies of $\vA$ and $\vA'$ coincide (for the
monodromies to exist, $\vA$ and $\vA'$ have to satisfy the assumptions
of Corollary \ref{co:infty}).
\end{remark}

The simplest class of modifications is the so-called \emph{elementary modifications}:

\begin{definition} Suppose the rational isomorphism $\vR:\vL\isorat\vL'$ is regular and has exactly one simple
zero. In this case, $\vA'$ is an \emph{elementary upper modification} of $\vA$, and $\vA$ is an \emph{elementary lower
modification} of $\vA'$.
\end{definition}

Note that an upper elementary modification $\vR:\vL\to\vL'$ is uniquely determined by the pair $(x,l)$, where $x\in\p1$
is the only zero of $\vR$ and the one-dimensional subspace $l\subset\vL_x$ is given by
$l=\ker(\vR(x):\vL_x\to\vL'_x)\subset\vL_x$.
Conversely, any pair $(x\in\p1,l\subset\vL_x)$ defines an elementary upper modification. Similarly, elementary lower
modifications of $\vL'$ are in on-to-one correspondence with pairs $(x,l')$ where $x\in\p1$, $l'\subset\vL'_x$ is a subspace
of codimension one (for $\vR:\vL\to\vL'$, $x$ is the only zero of $R$ and $l'=\im(\vR(x):\vL_x\to\vL'_x)$).

\begin{proposition} Suppose $(\vL,\vA)\in M_\theta$ for $\theta=(a_1,\dots,a_k;\rho_1,\rho_2,d_1,d_2;n)$,
and $\rho_1\ne\rho_2$. Let $(\vL',\vA')$ be an elementary upper
modification of $\vL$ given by $(x\in\p1;l\subset\vL_x)$. Then the
only cases when $(\vL',\vA')$ belongs to $M_{\theta'}$ for some
$\theta'\in\param$ are as follows:
\begin{enumerate}
\item If $x=\infty$, then $l$ must be an
eigenspace of $\vA_n:\vL_\infty\to\vL_\infty$ (the leading term of
$\vA=\vA_nz^n+\text{lower order terms}$). If, for instance,
$l=\ker(\vA_n-\rho_1)\subset L_{\infty}$, then
$\theta'=(a_1,\dots,a_k;\rho_1,\rho_2,d_1-1,d_2;n)$,
and an analogous formula holds when $l=\ker(\vA_n-\rho_2)$.

\item If $x=a_i$ is a zero of $\vA$ and $x-1\ne
a_j$ is not, then $l$ must be the kernel of
$\vA(x):\vL_x\to\vL_{x+1}$; in this case,
$\theta'=(a_1,\dots,a_i-1,\dots,a_k;\rho_1,\rho_2,d_1,d_2;n)$.
\end{enumerate}

In either case, the elementary modifications define an isomorphisms
$M_\theta\iso M_{\theta'}$.
\qed
\label{pp:modv}
\end{proposition}

\begin{remark} Sometimes an elementary modification
of a d-connection of type $\theta$ has simple poles,
which can be turned into simple zeroes using multiplication by a
scalar (for example, this happens if neither $x$ nor $x-1$ is a
pole). However, this procedure does not lead to an isomorphism between the
moduli spaces $M_\theta$ (at least assuming \eqref{cn:pole}--\eqref{cn:rho}
hold), because the corresponding spaces have different
dimensions.
\end{remark}

Thus, elementary modifications (upper or lower) allow to
identify $M_{\theta'}$ and $M_\theta$ if $\theta'$ is obtained
from $\theta$ by adding or subtracting $1$ to one of $a_i$'s or
$d_i$'s, provided certain conditions hold. Composing such
identifications, we get other isomorphisms between $M_\theta$ for
different $\theta\in\param_k$.

The situation is particularly simple if $\theta$ satisfies the
conditions \eqref{cn:mod}, \eqref{cn:rho}. Then $M_\theta$ and
$M_{\theta'}$ are naturally isomorphic if $\theta'$ is obtained
from $\theta$ by adding integers to $a_i$'s and $d_i$'s. In other
words, we have a natural action of the group
$G=(\Z)^k\times(\Z)^2$ on $\param_k$, and for any
$\theta\in\param_k$ satisfying \eqref{cn:mod},
\eqref{cn:rho} (in particular, for any $\theta\in\paramv_k$), we
get isomorphisms $M_\theta\to M_{g\theta}$ for all $g\in G$.

\subsection{Irreducibility of d-connections}
\label{sc:irreducibility}
Let $\vA(z)$ be a d-connection on a vector bundle $\vL$ on $\p1$.
Assume that $\vA(z)$ is non-degenerate at infinity in the sense that
\eqref{cn:dpv} holds. Denote by
$(\rho_1,\dots,\rho_m,d_1,\dots,d_m;n)$ the formal type of $\vA(z)$
at infinity.

For the morphism $\vA(z):\vL_z\to\vL_{z+1}$, its determinant is a map
$\det\vA(z):\bigwedge^m\vL_z\to\bigwedge^m\vL_{z+1}$; in other words,
$\det\vA(z)$ is a d-connection on the line bundle $\det
\vL\eqd\bigwedge^m\vL$. It is easy to see that $\det\vL$ has formal type
$(\rho_1\rho_2\cdots\rho_m,d_1+\dots+d_m;mn)$ at infinity. Let
$a_1,\dots,a_k\in\A1$ and $b_1,\dots,b_l\in\A1$ be zeroes and
poles (counted with multiplicity), respectively, of $\det\vA(z)$ on
$\A1$.

\begin{lemma} The collection $(a_1,\dots,a_k;b_1,\dots,b_l;\rho_1,\dots,\rho_m,d_1,\dots,d_m;n)$ satisfies the following
equalities:

\begin{align*}
mn&=k-l\\
\deg(\vL)&=-\sum_{i=1}^md_i-\sum_{i=1}^k a_i+\sum_{i=1}^l b_i.
\end{align*}
\qed
\label{lm:typecond}
\end{lemma}

\begin{corollary} Let $(\vL,\vA)$ be a d-connection of type
\begin{equation*}
\theta=(a_1,\dots,a_k;\rho_1,\rho_2,d_1,d_2;n)\in\param.
\end{equation*}
Then $k=2n$ and $\deg(\theta)=\deg(\vL)$ (see \eqref{cn:degree} for the
definition of $\deg(\theta)$); in particular, $\deg(\theta)$ is an
integer. \qed \label{co:degree}
\end{corollary}

\begin{lemma} Suppose $\theta\in\param$ satisfies \eqref{cn:irred}.
Then any $(\vL,\vA)\in M_{\theta}$ is irreducible: there is no rank 1
subbundle $\ell\subset\vL$ such that $\vA(\ell_z)\subset\ell_{z+1}$
for all $z$. \label{lm:irred}
\end{lemma}
\begin{proof} (Both the statement and its proof are completely analogous
to \cite[Proposition 1]{AL}.) Suppose $\ell\subset\vL$ is an invariant
subbundle of rank $1$, so that $\vA$ induces a d-connection
$\vA|_\ell$ on $\ell$. All zeroes of $\vA|_\ell$ belong to
$\{a_1,\dots,a_k\}$; besides, the formal type of $\vA|_\ell$ at
infinity is either $(\rho_1,d_1;n)$ or $(\rho_2,d_2;n)$. Now Lemma
\ref{lm:typecond} leads to a contradiction.
\end{proof}

\begin{corollary} Suppose $(\vL,\vA)\in M_\theta$ and suppose that $\theta\in\param_{2n}$ satisfies \eqref{cn:irred}. If
$\vL\simeq\vO(n_1)\oplus\vO(n_2)$, then $|n_1-n_2|\le n$.
\label{co:height}
\end{corollary}

\begin{proof}
Without loss of generality, we can assume $n_1\ge n_2$. Let
$\ell\subset\vL$ be a rank 1 subbundle of degree $n_1$. Since
$(\vL,\vA)$ is irreducible, $\ell$ is not $\vA$-invariant, and so the
rational map $\alpha:\ell\to\vL\to s^*\vL\to s^*(\vL/\ell)$ is not
identically zero. Notice that $\alpha$ can have at most a pole of
order $n$ at $\infty$ (and no other poles); thus,
$n_1=\deg(\ell)\le n+\deg(L/\ell)=n+n_2$.
\end{proof}

\section{Difference $PV$}

In this section, we study $M_\theta$ for
\begin{equation*}
\theta=(a_1,a_2,a_3,a_4;\rho_1,\rho_2,d_1,d_2;2)\in\paramv_4.
\end{equation*}
We assume $\deg(\theta)=-1$ (that is, $-d_1-d_2-\sum_{i=1}^4a_i=-1$).
Using modifications, we can make this assumption without loss of generality.

\subsection{$M_\theta$ as a quotient}
Let $(\vL,\vA)\in M_\theta$. By Corollary \ref{co:height}, $\vL$ is isomorphic
to $\vO\oplus\vO(-1)$. Let us choose an isomorphism $\vS:\vO\oplus\vO(-1)\iso\vL$;
then $\vA$ induces the d-connection $\vS(z+1)^{-1}\vA(z)\vS(z)$ of type $\theta$
on $\vO\oplus\vO(-1)$. Such a d-connection can be written as a matrix
\begin{equation}
\mA=\begin{bmatrix}a_{11}& a_{12}\\a_{21}& a_{22}\end{bmatrix},\quad
\begin{array}{c}a_{11},a_{22}\in\Gamma(\p1,\vO(2))\\a_{12}\in\Gamma(\p1,\vO(3))\\a_{21}\in\Gamma(\p1,\vO(1)).\end{array}
\label{eq:connmat}
\end{equation}

Of course, $\vS$ is not unique: it can be composed with an automorphism
of $\vO\oplus\vO(-1)$. Such an automorphism can be written as a matrix
\begin{equation}
\mR=\begin{bmatrix}r_{11}&r_{12}\\0&
r_{22}\end{bmatrix},\quad
\begin{array}{c}r_{11},r_{22}\in\C-\{0\}\\r_{12}\in\Gamma(\p1,\vO(1)).\end{array}
\label{eq:gaugemat}
\end{equation}
If we replace $\vS$ with $\vS\circ\mR$, then $\mA$ is replaced with its d-gauge transform
\begin{equation}
\mR(z+1)^{-1}\mA(z)\mR(z).
\label{eq:gauge2}
\end{equation}

\begin{lemma} Let $\vA$ be a d-connection on $\vO\oplus\vO(-1)$; its matrix $\mA$ is of the form \eqref{eq:connmat}.
We claim that $\vA$ is of type $\theta$ if and only if $\mA$ satisfies the following conditions:
\begin{gather}
\det(\mA)=(z-a_1)(z-a_2)(z-a_3)(z-a_4)\rho_1\rho_2 \label{cn:det}\\
a_{11}+a_{22}(1+z^{-1})=(\rho_1+\rho_2)z^2+(d_1\rho_1+d_2\rho_2)z+t(z^{-1}), \label{cn:tr}
\end{gather}
where $t(z^{-1})\in\C[[z^{-1}]]$ is a Taylor series in $z^{-1}$.
\label{lm:connmat}
\end{lemma}
\begin{proof}
$\vA$ is of type $\theta$ if and only if it satisfies the two conditions of Definition \ref{df:type}. Let us
reformulate the conditions in terms of $\mA$.

Definition \ref{df:type}(1) is equivalent to the condition that
\begin{equation}
\det(\mA)=c(z-a_1)(z-a_2)(z-a_3)(z-a_4)\quad\text{for some }c\in\C-\{0\}\label{cn:det1}
\end{equation}
(here we use that $\det(\mA)$ is a polynomial of degree $4$ in $z$). Now set
\begin{equation*}
\mS(z)\eqd\begin{bmatrix}1&0\\0&z^{-1}\end{bmatrix}
\end{equation*}
($\mS$ is a essentially a basis of $\vO\oplus\vO(-1)$ in a neighborhood of $\infty\in\p1$).
By Remark \ref{rm:characteristic polynomial}, Definition \ref{df:type}(2) is equivalent to the following two conditions:
\begin{align}
\det(\mS(z+1)^{-1}\mA(z)\mS(z))&=\rho_1\rho_2 z^4+\rho_1\rho_2(d_1+d_2)z^3+t_1(z^{-1})z^2\label{cn:det2}\\
\tr(\mS(z+1)^{-1}\mA(z)\mS(z))&=(\rho_1+\rho_2)z^2+(d_1\rho_1+d_2\rho_2)z+t_2(z^{-1})\label{cn:tr1}.
\end{align}
Here $t_1,t_2$ are Taylor series in $z^{-1}$.

It is easy to see that \eqref{cn:det} is equivalent to the combination of \eqref{cn:det1} and \eqref{cn:det2}
(here we use that
$\deg(\theta)=-1$), and \eqref{cn:tr} is equivalent to \eqref{cn:tr1}.
\end{proof}

\begin{corollary}
Denote by $X_\theta$ the space of matrices $\mA$ of the form \eqref{eq:connmat} that satisfy \eqref{cn:det} and
\eqref{cn:tr}; denote by $G$ be the group of matrices $\mR$ of the form \eqref{eq:gaugemat}. Let $G$ act on
$X_\theta$ via d-gauge transformations \eqref{eq:gauge2}. Then the quotient $X_\theta/G$ is canonically
isomorphic to $M_\theta$.\qed
\label{co:dPV as quotient}
\end{corollary}

\subsection{Geometric description of $M_\theta$}
\label{sc:Geometry of dPV}
In this section, we will derive Theorem \ref{th:Geometry of dPV} from another geometric description of $M_\theta$
(Theorem \ref{th:Geometry of dPV tilde}).
Recall that Theorem \ref{th:Geometry of dPV} realizes $M_\theta$ as an open subset of a blow-up of $(\p1)^2$;
in Theorem \ref{th:Geometry of dPV tilde}, we use a different rational surface in place of $(\p1)^2$. Of the two descriptions,
Theorem \ref{th:Geometry of dPV tilde} uses somewhat more natural constructions (however, see Remark \ref{rm:dPV and Gamma});
for instance, all four points $a_1,\dots,a_4$ appear in a symmetric manner. On the other hand, the advantage of Theorem \ref{th:Geometry of dPV}
is that $(\p1)^2$ has natural coordinates $(q,p)$, which can then be viewed as `rational coordinates' $q,p:M_\theta\to\p1$.
This makes Theorem \ref{th:Geometry of dPV} more suitable for writing formulas.

As before, $(\vL,\vA)\in M_\theta$, $\vS:\vO\oplus\vO(-1)\iso\vL$, and $\mA$ is the matrix of $\vA$ relative
to $\vS$.
Notice that the matrix element $a_{21}\in\Gamma(\p1,\vO(1))$ is not
identically zero, because $(\vL,\vA)$ is irreducible. Therefore,
$a_{21}$ has a single zero on $\p1$; let us denote it by $q\in\p1$.
Set $\tilde p\eqd a_{11}(q)\in (\vO(2))_q$.

\begin{proposition} $\tilde p$ and $q$ depend only on $(\vL,\vA)\in M_\theta$ and not on $\vS$.
\end{proposition}
\begin{proof} This statement can be easily checked directly by calculating the d-gauge
transformation \eqref{eq:gauge2} with the gauge matrix \eqref{eq:gaugemat}.
It is also possible to provide a geometric explanation in the
spirit of \cite[Section 4.1]{AL}.
\end{proof}

$(q,\tilde p)$ can be viewed as a map $\tilde P:M_\theta\to
\tK$, where $\tK\eqd\V(\vO(2)^\vee)$ is the total
space of the line bundle $\vO(2)$. As we will see, the map $\tilde
P:M_\theta\to\tK$ is a regular birational morphism. Since $M_\theta$ is
a smooth algebraic surface, $\tP$ identifies $M_\theta$
with an open subset of a blow-up of $\tK$. Let us describe the
blow-up.

Let us start with some general remarks about the geometry of
$\tK$. Clearly, $\tK$ is fibered over $\p1$ so that the fiber over
$z\in\p1$ is $\vO(2)_z$. If $f$ is a (rational) section of $O(2)$
that is regular at $z$, then its value $f(z)\in\vO(2)_z$ can be
viewed as a point of $\tK$; we will denote this point by
$(z,f(z))$. For example, $(z,0(z))$ is the zero element in the
fiber of $\tK$ over $z\in\p1$.

Now let $\tilde\sigma_c:\tK_c\to\tK$ be the blow-up of $\tK$ at
$c\eqd(z,f(z))$. Then the exceptional divisor
$\tilde\sigma_c^{-1}(c)\subset\tK_c$ is isomorphic to the
projective line $\PP(T_c\tK)$; that is, points of
$\tilde\sigma_c^{-1}(c)$ correspond to lines in the tangent space to
$\tK$ at $c$. Any smooth curve $C\subset\tK$ that passes through
$c$ defines such a line (the tangent line to $C$ at $c$). In
particular, we can take $C$ to be the graph $\{(x,f(x)):x\in\p1\}$
of $f$; denote the corresponding point of $\tK_c$ by $(z,f'(z))$.
Any other rational section $g$ of $\vO(2)$ defines a point
$(z,g'(z))\in\tK_c$ provided $g$ is regular at $z$ and
$g(z)=f(z)$.

\begin{theorem}
\begin{enumerate}
\item The map $\tilde P:M_\theta\to\tK$ is a regular birational
morphism of smooth algebraic surfaces.

\item Let $\tilde\sigma_1:\tK_1\to\tK$ be the blow-up of $\tK$ at
the following $6$ points: $(a_i,0(a_i))$ ($i=1,\dots,4$) and
$(\infty,(\rho_jz^2)(\infty))$ ($j=1,2$). Let $\sigma_2:\tK_2\to
\tK_1$ be the blow-up of $\tK_1$ at the two points
$(\infty,(\rho_jz^2+\rho_j d_jz)'(\infty))$, $j=1,2$ (these points
belong to the preimages of $(\infty,(\rho_jz^2)(\infty))$,
$j=1,2$). Then the map $\tP$ induces an open embedding
$\tP_2:M_\theta\hookrightarrow \tK_2$.

\item The complement to $\tP_2(M_\theta)$ in $\tK_2$ is the union
of the proper preimages of the following curves: the zero section
$\{(z,0(z)):z\in\p1\}\subset\tK$, the fiber at infinity
$\{(\infty,az^2(\infty)):a\in\C\}\subset\tK$, and two
exceptional curves
$\tilde\sigma_1^{-1}(\infty,(\rho_jz^2)(\infty))\subset\tK_1$.
\end{enumerate}
\label{th:Geometry of dPV tilde}
\end{theorem}

The proof of Theorem \ref{th:Geometry of dPV tilde} is given in Section \ref{sc:Geometry of dPV tilde}.
Let us now derive Theorem \ref{th:Geometry of dPV} from Theorem \ref{th:Geometry of dPV tilde}.

\begin{proof}[Proof of Theorem \ref{th:Geometry of dPV}]
For $(\vL,\vA)\in M_\theta$, consider the expression
\begin{equation}
p\eqd\dfrac{\tilde p}{(q-a_3)(q-a_4)}.
\label{eq:p dPV}
\end{equation}
Here the denominator is the value of the section
$(z-a_3)(z-a_4)\in\Gamma(\p1,\vO(2))$ at $z=q\in\p1$. Both the numerator and the denominator are elements
of $\vO(2)_q$; therefore, $p\in\C$ provided the denominator does not vanish. We can view $p$ as
a rational mapping $p:M_\theta\to\p1$. Actually, Theorem \ref{th:Geometry of dPV tilde} implies
that $p:M_\theta\to\p1$ is regular: the corresponding rational mapping $\tK\dashrightarrow\p1$ has singularities at
$(a_3,0(a_3)), (a_4,0(a_4))\in\tK$, but the blow-up $\tK_1\to\tK$ resolves the singularities. We therefore obtain a regular
mapping $P\eqd(q,p):M_\theta\to(\p1)^2$. We claim that $P$ induces an embedding $P_2:M_\theta\hookrightarrow K_2$,
where $K_2$ is the blow-up of $(\p1)^2$ described in Theorem \ref{th:Geometry of dPV}.

Let us consider the birational mapping $\Phi:(q,\tilde p)\mapsto(q,p):\tK\dashrightarrow(\p1)^2$. It is easy to see that
$\Phi$ induces an open embedding $\Phi_1:\tK_1\hookrightarrow K_1$, and the complement $K_1-\Phi(\tK_1)$ is the proper
preimage of $\p1\times\{\infty\}\subset(\p1)^2$ under the blow-up $K_1\to(\p1)^2$. To complete the proof, we should now check
that $\Phi_1$ maps the centers of the blow-up $\tK_2\to\tK_1$
to the centers of the blow-up $K_2\to K_1$. This also follows from the formulas.
\end{proof}

\begin{remark}\label{rm:dPV and Gamma}
Geometrically, formula \eqref{eq:p dPV} can be explained as a multiplication a d-connection by a scalar.
For $(\vL,\vA)\in M_\theta$, consider the d-connection
\begin{equation*}
\widetilde\vA\eqd\dfrac{1}{(z-a_3)(z-a_4)}\vA
\end{equation*}
on $\vL$. Then $\widetilde\vA$ has simple zeroes at $a_1,a_2$, simple poles at $a_3,a_4$, and its formal type at infinity
is $(\rho_1,\rho_2;d_1+a_3+a_4,d_2+a_3+a_4;0)$. Moreover, we can then view $M_\theta$ as the moduli space of d-connections
of this kind (as in Remark \ref{rm:poles and zeroes}). For d-connections of this kind,
$p$ plays the role of $\tilde p$, and Theorem \ref{th:Geometry of dPV}
plays the role of Theorem \ref{th:Geometry of dPV tilde}.
\end{remark}

\subsection{Proof of Theorem \ref{th:Geometry of dPV tilde}}\label{sc:Geometry of dPV tilde}
The most direct way to prove Theorem \ref{th:Geometry of dPV tilde} is by bringing matrices \eqref{eq:connmat} to
some `normal form'. We will not reproduce all calculations here; the idea of the proof is as follows:

Denote by $\tM_\theta$ the open subset of $\tK_2$ described in Theorem \ref{th:Geometry of dPV tilde}(3) (that is,
the complement of proper preimages of the zero section, the fiber at infinity, and two exceptional curves). We need
to show that the map $\tP:M_\theta\to\tK$ lifts to an isomorphism $M_\theta\to\tM_\theta$. Let us consider open
sets
\begin{align*}
U_0&\eqd q^{-1}(\p1-\{\infty\})\subset M_\theta &
U_\infty&\eqd q^{-1}(\p1-\{0\})\subset M_\theta\\
\tU_0&\eqd q^{-1}(\p1-\{\infty\})\subset\tM_\theta &
\tU_\infty&\eqd q^{-1}(\p1-\{\infty\})\subset\tM_\theta.
\end{align*}
It suffices to show that $\tP$ lifts to isomorphisms $U_0\iso\tU_0$, $U_\infty\iso\tU_\infty$. We will show this
by writing $U_0$, $U_\infty$ explicitly as zero loci of polynomial equations.

Let $(\vL,\vA)$ be a point of $U_0$. Then $q=q(\vL,\vA)\in\C$ and $\tilde
p=\tilde p(\vL,\vA)\in(\vO(2))_q=\C$. It is easy to see that there
exists an isomorphism $\vS:\vO\oplus\vO(-1)\iso\vL$, unique up to
a multiplicative constant, such that the matrix of the
d-connection $\vA$ relative to $\vS$ is
\begin{equation}
\mA=\begin{bmatrix}a_{11}=\tilde p& a_{12}\\a_{21}=z-q&
a_{22}\end{bmatrix},\quad
\begin{array}{c}a_{22}\in\Gamma(\p1,\vO(2))\\a_{12}\in\Gamma(\p1,\vO(3)).\end{array}
\label{eq:connmat0}
\end{equation}
Essentially, \eqref{eq:connmat0} serves as a normal form of
d-connections $(\vL,\vA)$ (provided $q\ne\infty$). The conditions
\eqref{cn:det}, \eqref{cn:tr} now become equations on $a_{12}$, $a_{22}$.
Explicitly, $a_{12}$ and $a_{22}$ are determined by
their coefficients
\begin{align*}
a_{12}&=a_{12,3}z^3+a_{12,2}z^2+a_{12,1}z+a_{12,0}\\
a_{22}&=a_{22,2}z^2+a_{22,1}z+a_{22,0},
\end{align*}
 and \eqref{cn:det}, \eqref{cn:tr} is a system of polynomial equations on
$a_{i2,j}$, $\tilde p$, and $q$. Solving these equations, we find
polynomial (in $\tilde p$ and $q$) formulas for all $a_{i2,j}$,
except for $r\eqd a_{22,0}$. The equation on $r$ looks as follows:
\begin{equation}
\tilde pr=F(\tilde p,q), \label{eq:pqr}
\end{equation}
where $F(\tilde p,q)$ is a polynomial. Thus, $U_0$ is
identified with the zero locus of the equation \eqref{eq:pqr} in
the three-dimensional space with coordinates $\tilde p,q,r$.

Besides, $F(0,q)=c(q-a_1)(q-a_2)(q-a_3)(q-a_4)$ for some $c\in\C-\{0\}$.
Therefore, the map $(\tilde p,q):U_0\to\A2$ identifies $U_0$ with
the complement to the proper preimage of `the $q$-axis'
$\{(0,q)\}$ in the blow-up of $\A2$ at the four points $(\tilde
p,q)=(0,a_i)$, $i=1,\dots,4$. This complement is exactly $\tU_0$.

Similar approach works for $U_\infty$. For $(\vL,\vA)\in
U_\infty$, set $\omega\eqd(q(\vL,\vA))^{-1}\in\C$, $\pi\eqd\frac{\tilde
p(\vL,\vA)}{q(\vL,\vA)^2}\in\C$, where the denominator is understood as
the value of $z^2\in\Gamma(\p1,\vO(2))$ at $z=q$. One can think of $\omega$ and $\pi$
as the coordinates on the complement to the zero locus of $q$ in
$\tK$. Then there is a unique up to a multiplicative constant
choice of $\vS:\vO\oplus\vO(-1)\iso\vL$ such that the matrix of $\vA$
is
\begin{equation*}
\mA=\begin{bmatrix}\pi z^2&a_{12}\\
1-\omega z&a_{22}\end{bmatrix}.
\end{equation*}
Again, we get a system of polynomial equations on the coefficients
of $a_{i2}$. Solving the equations, we find polynomial (in $\pi$
and $\omega$) formulas for all $a_{i2,j}$, except for
$r= a_{22,0}$. In this case, the equation on $r$ is
\begin{equation}
\pi\omega^2r=G(\pi,\omega), \label{eq:piomegar}
\end{equation}
where $G(\pi,\omega)$ is a polynomial. Therefore, $U_\infty$ is
the zero locus of the equation \eqref{eq:piomegar} in the
three-dimensional space with coordinates $\pi,\omega,r$. Again, from the
formula for $G(\pi,\omega)$, one easily sees the isomorphism $U_\infty\iso\tU_\infty$.

For instance, let us consider the neighborhood of $\omega=0$ (the complement of $\omega=0$ is covered
by $U_0$). One can check that $G(\pi,0)=(\pi-\rho_1)(\pi-\rho_2)$, so when $\omega=0$, either $\pi=\rho_1$,
or $\pi=\rho_2$. Consider the neighborhood of the set $\omega=0$, $\pi=\rho_1$ in $U_\infty$.
It follows that $\pi_1\eqd(\pi-\rho_1)/\omega$ is a regular function on the neighborhood
($\pi_1$ is the coordinate on the blow-up of the $\omega$-$\pi$ plane at $(\omega,\pi)=(0,\rho_1)$).
We can then rewrite \eqref{eq:piomegar} in variables $\pi_1$, $\omega$, and $r$:
\begin{equation*}
(\omega\pi_1+\rho_1)r\omega=H(\pi_1,\omega),
\end{equation*}
where $H(\pi_1,\omega)$ is a polynomial such that $H(\pi_1,0)=(\rho_2-\rho_1)(\pi_1-\rho_1 d_1)$; therefore,
$r$ is essentially the coordinate on the blow-up of the $\omega$-$\pi_1$ plane at $(\omega,\pi_1)=(0,\rho_1 d_1)$.
Of course, the neighborhood of the set $\omega=0$, $\pi=\rho_2$ in $U_\infty$ has a similar description.
\qed

\begin{remark} Theorem \ref{th:Geometry of dPV} can be also proved in a more geometric way,
in the spirit of \cite[Theorem 3]{AL}.
\end{remark}

\subsection{Proof of Theorem \ref{th:dPV}}\label{sc:dPV}
The proof of Theorem \ref{th:dPV} is also based on calculations. The calculations are simplified
by the observation that it suffices to check the formulas \eqref{eq:dPV} on a dense subset of $M_\theta$; we can therefore assume that
$q,q'\ne\infty$.

Take $(\vL,\vA)\in M_\theta$ and set $(\vL',\vA')\eqd dPV(\vL,\vA)$. Let us assume $q(\vL,\vA)\ne\infty$ (that is, $(\vL,\vA)\in U_0$),
then there is an isomorphism $\vS:\vO\oplus\vO(-1)\iso\vL$ such that the matrix of $\vA$ relative to $\vS$ is of the form \eqref{eq:connmat0}.
Using the formula $\tilde p=p(q-a_3)(q-a_4)$, we can write the matrix as
\begin{equation}
\mA=\begin{bmatrix}p(q-a_3)(q-a_4)& a_{12}\\z-q&a_{22}\end{bmatrix},\quad
\begin{array}{c}a_{22}\in\Gamma(\p1,\vO(2))\\a_{12}\in\Gamma(\p1,\vO(3)).\end{array}
\label{eq:connmatV}
\end{equation}
Recall also that $a_{12}$, $a_{22}$ are polynomials of $z$ whose coefficients are rational functions of $p$, $q$.

Similarly, if we assume $q(\vL',\vA')\ne0$, there exists an isomorphism $\vS':\vO\oplus\vO(-1)\iso\vL'$ such that the matrix of $\vA'$ relative to
$\vS'$ is of the form
\begin{equation}
\mA'=\begin{bmatrix}p'(q'-a_3)(q'-a_4)& a'_{12}\\z-q'&a'_{22}\end{bmatrix},\quad
\begin{array}{c}a'_{22}\in\Gamma(\p1,\vO(2))\\a'_{12}\in\Gamma(\p1,\vO(3)).\end{array}
\label{eq:connmatV'}
\end{equation}
By the definition of $dPV$, the matrix $\mA'$ is the d-gauge transformation of $\mA$:
\begin{equation}
\mA'(z)=\mR(z+1)^{-1}\mA(z)\mR(z),
\label{eq:gaugeV}
\end{equation}
where $\mR$ is the matrix of the rational map $\vR:\vL'\isorat\vL$ (from the definition of $dPV$)
with respect to the bases $\vS$, $\vS'$. It follows from the properties of modifications (Section \ref{sc:operations})
that $\vR$ induces a regular map $\vL'\to\vL\otimes\vO(1)$ whose determinant has simple zeroes at $a_1$, $a_2$ and
no other zeroes. In other words, $\mR$ is of the form
\begin{equation*}
\mR=\begin{bmatrix}r_{11}&r_{12}\\r_{21}&
r_{22}\end{bmatrix},\quad
\begin{array}{c}r_{11},r_{22}\in\Gamma(\p1,\vO(1))\\r_{21}\in\C,r_{12}\in\Gamma(\p1,\vO(2)),\end{array}
\label{eq:gaugematV}
\end{equation*}
such that
\begin{equation}
\det(\mR)=c(z-a_1)(z-a_2)\quad(c\in\C-\{0\}).
\label{eq:gaugematdetV}
\end{equation}

\eqref{eq:gaugematdetV} yields polynomial equations on the coefficients of $r_{11},r_{12},r_{21},r_{22}$;
the condition that \eqref{eq:gaugeV} gives a matrix $\mA'$ of the form \eqref{eq:connmatV'} also gives such equations. The resulting system
determines $\mR$ up to a multiplicative constant. From \eqref{eq:gaugeV}, we now obtain a formula for
the matrix $\mA'$ in terms of $p$ and $q$; in particular, we  can derive \eqref{eq:dPV}.\qed

\section{Difference $PV$ and classical $PVI$}

\subsection{Geometry of $PVI$}
\label{sc:Geometry of PVI}
Let us recall the description of the surface $M_{(x,\lambda)}$. We will suppose that
\begin{equation}
\sum_{i=1}^4(\lambda_i^-+\lambda_i^+)=1.
\label{eq:lambdas}
\end{equation}
It is easy to see that $M_{(x,\lambda)}$ only depends on the classes of $\lambda_i^\pm$
in $\C/\Z$ (because of modifications of bundles with connections), so our assumption
does not restrict the generality.

Suppose $x\in X$, $\lambda\in\Lambda$, and let $K_x$ be the total space of the line bundle
$\Omega_\p1(x_1+\dots+x_4)$. Let $b_i\subset K_x$ be the fiber over $x_i\in\p1$. Notice that
the residue of 1-forms identifies the fiber of $\Omega(x_1+\dots+x_4)$ over $x_i$ with $\C$,
so we get a canonical isomorphism $\res_i:b_i\iso\A1$. Denote by $\tilde M_{(x,\lambda)}$ the blow-up
of $K_x$ at the eight points $(\res_i)^{-1}(\lambda_i^\pm)$, $i=1,\dots,4$, and let
$M'_{(x,\lambda)}\subset \tilde M_{(x,\lambda)}$ be the complement to the proper preimages of $b_i\subset K_x$.

\begin{proposition}
There exists an isomorphism $M_{(x,\lambda)}\iso M'_{(x,\lambda)}$. \qed
\label{pp:Geometry of PVI}
\end{proposition}

Proposition \ref{pp:Geometry of PVI} is a slight generalization of \cite[Theorem 3]{AL} (see also \cite[Theorem 2.2]{IIS2}):
\cite{AL} works only with ${\mathrm SL}(2)$-bundles, which corresponds
to assuming $\lambda_i^-+\lambda_i^+=0$, ($i=2,3,4$). However, the general case is easily reduced to this special case. Let us sketch the
construction of the map $M_{(x,\lambda)}\iso M'_{(x,\lambda)}$.

Given $(\vL,\nabla)\in M_{(x,\lambda)}$, one can show that $\vL\simeq\vO\oplus\vO(-1)$ (this is similar to Corollary \ref{co:height}). If
we fix an isomorphism $\vO\oplus\vO(-1)\iso\vL$, the connection $\nabla$ is determined by its matrix
\begin{equation*}
\mM(z)=\begin{bmatrix}m_{11}&m_{12}\\m_{21}&m_{22}\end{bmatrix}\quad
\begin{array}{c}m_{11},m_{22}\in\Gamma(\p1,\Omega_\p1(x_1+\dots+x_4))\\m_{12}\in\Gamma(\p1,\Omega_\p1(x_1+\dots+x_4)\otimes\vO(1))
\\m_{21}\in\Gamma(\p1,\Omega_\p1(x_1+\dots+x_4)\otimes\vO(-1)).\end{array}
\end{equation*}
It can be proved that $m_{21}$ is not identically zero (because $(\vL,\nabla)$ is irreducible; this is similar to Lemma \ref{lm:irred}).
Therefore, $m_{21}$ has a single zero on $\p1$; denote it by $q^{PVI}$. Set $p^{PVI}\eqd m_{11}(q^{PVI})$. Note that
$p^{PVI}$ belongs to the fiber of $\Omega_\p1(x_1+\dots+x_4)$ over $q^{PVI}\in\p1$. In other words, $p^{PVI}$ is a point of the total
space $K_x$ (in the notation of Section \ref{sc:Geometry of dPV}, the point is $(q^{PVI},p^{PVI})\in K_x$). One can check that
$q^{PVI}$ and $p^{PVI}$ depend only on $(\vL,\nabla)$, not on the choice of $\vO\oplus\vO(-1)\iso\vL$. Therefore,
we obtain a regular
map $M_{(x,\lambda)}\to K_x$. Proposition \ref{pp:Geometry of PVI} claims the map induces an isomorphism $M_{(x,\lambda)}\iso M'_{(x,\lambda)}$.

\begin{proof}[Proof of Theorem \ref{th:dPV and PVI}] Let $\theta\in\paramv_4$, $x\in X$, and $\lambda\in\Lambda$ be as in Theorem \ref{th:dPV and PVI};
we will define the isomorphism $M_\theta\to M_{(x,\lambda)}$ by explicit formulas. Let $q,p:M_\theta\to\p1$ be the `coordinates' from Theorem
\ref{th:Geometry of dPV}. Consider the expression
\begin{equation*}
p^{PVI}\eqd(z^{-1}dz)_{z=p}q,
\end{equation*}
where $(z^{-1}dz)_{z=p}\in(\Omega_\p1(x_1+\dots+x_4))_p$ is the value of $z^{-1}dz\in\Gamma(\p1,\Omega_\p1(x_1+\dots+x_4))$ at $z=p$. Then
$p^{PVI}\in(\Omega_\p1(x_1+\dots+x_4))_p$, provided $q\ne\infty$. Let us also set $q^{PVI}\eqd p$.

If $q\ne\infty$, we have $(q^{PVI},p^{PVI})\in K_x$; in this manner, we get a rational map
\begin{equation*}
M_\theta\dashrightarrow K_x:(q,p)\mapsto(q^{PVI},p^{PVI}).
\end{equation*}
Using Theorem \ref{th:Geometry of dPV} and Proposition \ref{pp:Geometry of PVI}, it is easy to see that the map is actually regular, and that it
lifts to an isomorphism $M_\theta\to M_{(x,\lambda)}$.
\end{proof}

\subsection{Classical $PVI$}
The isomonodromy deformation of bundles with connections gives a system of differential equations on the
`coordinates' $q^{PVI}$, $p^{PVI}$ (the `usual' $PVI$). Here $q^{PVI}$, $p^{PVI}$ are
viewed as functions of $x_1,\dots,x_4$, while $\lambda_i^\pm$ are fixed parameters. Let us recall the explicit
formulas (which we adapted from \cite{IIS}).

For simplicity, we assume, in addition to \eqref{eq:lambdas}, that $x_4=\infty$. Define the new parameters
by $\kappa_i\eqd\lambda_i^+-\lambda_i^-$, $i=1,\dots,4$, and let us replace the variable $p^{PVI}$ with
\begin{equation*}
\tilde p^{PVI}\eqd(p^{PVI}/dz)-\sum_{i=1}^3\dfrac{\lambda_i^-}{z-x_i}.
\end{equation*}
Since $p^{PVI}\in(\Omega_\p1(x_1+\dots+x_4))_{q^{PVI}}$, the ratio $p^{PVI}/dz$ (if it is defined) is a number. The
advantage of $\tilde p^{PVI}$ is that the differential equations for $q^{PVI}$, $\tilde p^{PVI}$ involve
fewer parameters: $\kappa_i$'s, rather than $\lambda_i^\pm$'s.

Set also
\begin{equation*}
\kappa_0\eqd\frac{1}{2}\left(1-\sum_{i=1}^4\kappa_i\right),
\end{equation*}
and $q_i\eqd q^{PVI}-x_i$, $i=1,2,3$.
Define the Hamiltonians $h_i$, $i=1,2,3$ by
\begin{equation*}
h_i\eqd\frac{(q_1q_2q_3)(\tilde p^{PVI})^2-\left((\kappa_i-1)q_jq_k+\kappa_jq_iq_k+\kappa_kq_iq_j\right)\tilde p^{PVI}+
\kappa_0(\kappa_0+\kappa_4)}{(x_i-x_j)(x_i-x_k)}.
\end{equation*}
The equations can then be written in the Hamiltonian form as
\begin{equation}
\frac{\partial q^{PVI}}{\partial x_i}=\frac{\partial h_i}{\partial \tilde p^{PVI}},\qquad
\frac{\partial \tilde p^{PVI}}{\partial x_i}=\frac{\partial h_i}{\partial q^{PVI}},\qquad (i=1,2,3).
\label{eq:standardPVI}
\end{equation}

The system \ref{eq:standardPVI} can be reduced to the usual form of $PVI$ as follows: set
\begin{equation*}
y\eqd\frac{q^{PVI}-x_1}{x_2-x_1},\qquad x\eqd\frac{x_3-x_1}{x_2-x_1}.
\end{equation*}
Then \ref{eq:standardPVI} implies that $y$ depends only on $x$, not on $x_1$, $x_2$, $x_3$, and that
$y$ satisfies the $PVI$ equation
\begin{multline}
\dfrac{d^2y}{dx^2}=\dfrac{1}{2}\left(\dfrac{1}{y}+\dfrac{1}{y-1}+\dfrac{1}{y-x}\right)\left(\dfrac{dy}{dx}\right)^2-
\left(\dfrac{1}{x}+\dfrac{1}{x-1}+\dfrac{1}{y-x}\right)\dfrac{dy}{dx}+
\cr
\dfrac{y(y-1)(y-x)}{x^2(x-1)^2}\left(\kappa_4^2-\kappa_1^2\dfrac{x}{y^2}+\kappa_2^2\dfrac{x-1}{(y-1)^2}+
(1-\kappa_3^2)\dfrac{x(x-1)}{(y-x)^2}\right).
\end{multline}

\subsection{Isomonodromy deformation of d-connections}
\label{sc:isomonodromic}
Let us prove Theorem \ref{th:isomonodromic}(1).
Informally, we need to show that, given $\theta\in\paramv_{2n}$ and $\rho'_1,\rho'_2\in\C$,
any d-connection of type $\theta$ has a natural `first order deformation' that is `of type'
\begin{equation*}
\theta^\epsilon\eqd(a_1,\dots,a_{2n};\rho_1+\epsilon\rho'_1,\rho_2+\epsilon\rho'_2,d_1,d_2;n).
\end{equation*}
Here $\epsilon$ is the parameter of the deformation, and all calculations are done
modulo $\epsilon^2$, that is, over the ring of dual numbers
$\C^\epsilon\eqd\C[\epsilon]/(\epsilon^2)$.
First, let us prove a `formal' statement:

\begin{proposition} \label{pp:isomonodromic:formal}
Suppose the matrix $\mA(z)=\sum_{i\le n}\mA_iz^i$ over $\C((z^{-1}))$ has formal type
$(\rho_1,\dots,\rho_m;d_1,\dots,d_m;n)$ at infinity (see Proposition \ref{pp:infty}).
For any collection $\rho_1',\dots,\rho_m'\in\C$, there exists a
gauge matrix $\mR^\epsilon(z)=\mR(z)+\epsilon\mR'(z)$, where $\mR(z)$ is as in Proposition \ref{pp:infty}
(that is, $\mR(z)$ is an invertible $m\times m$ matrix over $\C[[z^{-1}]]$), and $\mR'(z)$ is an
$m\times m$ matrix over the ring of formal Laurent series $\C((z^{-1}))$
such that
\begin{equation}
\mR^\epsilon(z+1)^{-1}\mA(z)\mR^\epsilon(z)=\diag((\rho_1+\rho'_1\epsilon)(z^n+d_1z^{n-1}),\dots,
(\rho_m+\epsilon\rho'_m)(z^n+d_m z^{n-1})).
\label{cn:formaltypee}
\end{equation}
The matrix $\mR^\epsilon(z)$ is unique up to right multiplication by a diagonal matrix with entries in $\C^\epsilon$.
\end{proposition}
\begin{proof}
\eqref{cn:formaltypee} is equivalent to the following two conditions:
\begin{gather}
\mR(z+1)^{-1}\mA(z)\mR(z)=\diag(\rho_1z^n+\rho_1d_1z^{n-1},\dots,
\rho_mz^n+\rho_md_mz^{n-1})\label{eq:power0}\\
\begin{gathered}
\mR(z+1)^{-1}\mA(z)\mR'(z)-\mR(z+1)^{-1}\mR'(z+1)\mR(z+1)^{-1}\mA(z)\mR(z)=\\
\diag(\rho'_1z^n+\rho'_1d_1z^{n-1},\dots,\rho'_mz^n+\rho'_md_mz^{n-1}).
\end{gathered}\label{eq:power1}
\end{gather}
As $\mA(z)$ has formal type $(\rho_1,\dots,\rho_m;d_1,\dots,d_m;n)$ at infinity,
there exists a matrix $\mR(z)$ satisfying \eqref{eq:power0}; moreover, $\mR(z)$
is unique up to right multiplication by a constant diagonal matrix (Proposition \ref{pp:infty}).
Once \eqref{eq:power0} is satisfied, \eqref{eq:power1} can be rewritten as
\begin{equation}
\label{eq:power1'}
\mB(z)\mS(z)-\mS(z+1)\mB(z)=\diag(\rho'_1z^n+\rho'_1d_1z^{n-1},\dots,\rho'_mz^n+\rho'_md_mz^{n-1}),
\end{equation}
where we set $\mB(z)\eqd\diag(\rho_1z^n+\rho_1d_1z^{n-1},\dots,
\rho_mz^n+\rho_md_mz^{n-1})$, and $\mS(z)\eqd\mR(z)^{-1}\mR'(z)$. One can view
\eqref{eq:power1'} as a difference equation on the matrix $\mS(z)$; it is easy to see
that the only solutions whose matrix elements are Laurent series are given by
$\mS(z)=\diag((\rho'_1/\rho_1)z+c_1,\dots,(\rho'_m/\rho_m)z+c_m)$, where $c_i$'s are arbitrary constants.
This implies the statement.
\end{proof}

Proposition \ref{pp:isomonodromic:formal} allows us to construct the natural first order deformation,
thus proving Theorem \ref{th:isomonodromic}(1). The construction is most easily described using the following well-known
statement.
\begin{lemma}
Let $\vL$ be a vector bundle on $\p1$ and let $\vS(z):\C^2\iso\vL_z$ be a trivialization of $\vL$ in the punctured
formal neighborhood of $\infty$ (so $\vS(z)$ is essentially a matrix whose entries belong to $\C((z^{-1}))$). Then
there exists a unique vector bundle $\vL^\vS$ such that $\vL$ and $\vL^\vS$ have equal restrictions to $\p1-\{\infty\}$
and that the map
\begin{equation*}
\vS(z):\C^2\to\vL_z=(\vL^\vS)_z
\end{equation*}
extends to a trivialization of $\vL^\vS$ in the formal neighborhood of $\infty$. \qed
\label{lm:mod via matrices}
\end{lemma}
Notice that
Lemma \ref{lm:mod via matrices} still works when $\vS$ depends on parameters. In this case,
the modification $\vL^\vS$ will also depend on the parameters.

\begin{proof}[Proof of Theorem \ref{th:isomonodromic}(1)]
Take $\rho=(\rho_1,\rho_2)\in P$,
$(\vL,\vA)\in M_{\theta(\rho)}$. Take a tangent vector
$\tau=\rho'_1\frac{\partial}{\partial\rho_1}+\rho'_2\frac{\partial}{\partial\rho_2}$
to $P$ at $\rho$. Let us construct a natural lifting of $\tau$ to
a tangent vector a tangent vector $\tau_M$ to $M$ at $(\vL,\vA)\in M$.

Choose a trivialization $\vS(z):\C^2\iso\vL_z$ on the
neighborhood of $\infty\in\p1$. The matrix
\begin{equation*}
\mA(z)\eqd\vS^{-1}(z+1)\vA(z)\vS(z)
\end{equation*}
of $\vA$ relative to $\vS$ satisfies the assumption of Proposition \ref{pp:isomonodromic:formal}.
Let us set $\vS^\epsilon(z)\eqd\vS(z)\mR^\epsilon(z)$, where the matrix $\mR^\epsilon(z)$ is given by Proposition
\ref{pp:isomonodromic:formal}. We can view $\vS^\epsilon(z)$ as a trivialization of $\vL$ in the punctured formal
neighborhood of $\infty\in\p1$ that depends on $\epsilon\in\C^\epsilon$. Lemma \ref{lm:mod via matrices}
defines a vector bundle $\vL^\epsilon\eqd\vL^{\vS^\epsilon}$ that depends on $\epsilon$.

$\vL^\epsilon$ and $\vL$ coincide on $\p1-\{\infty\}$ (for any value of the parameter $\epsilon$), thus the d-connection
$\vA$ on $\vL$ induces a d-connection $\vA^\epsilon$ on $\vL^\epsilon$. Notice also that when $\epsilon=0$,
we have $\vL^\epsilon=\vL$, $\vA^\epsilon=\vA$. The pair $(\vL^\epsilon,\vA^\epsilon)$ define a tangent vector
$\tau_M$ to $M$ at $(\vL,\vA)$. The vector $\tau_M$ does not depend on the choice of $\mR^\epsilon$.
It is easy to see that as $\tau$ and $(\vL,\vA)$ vary, the lifting $\tau_M$
defines a flat algebraic connection on $M\to P$.
\end{proof}

\subsection{Isomonodromy deformation for $2n=4$}
Suppose now that $2n=4$, $\deg(\theta)=-1$. Then the construction of the previous section can be reformulated more
explicitly. Instead of working with d-connections, let us consider their
matrices (that is, we think of $M_\theta$ as a quotient $X_\theta/G$, see
Corollary \ref{co:dPV as quotient}).

Let $(\vL,\vA)$ and $(\vL^\epsilon,\vA^\epsilon)$ be as above. Choose a trivialization
$\vS^\epsilon:\vO\oplus\vO(-1)\iso\vL^\epsilon$ (depending on $\epsilon$). When $\epsilon=0$,
$\vS^\epsilon$ becomes a trivialization $\vS:\vO\oplus\vO(-1)\iso\vL$. Let $\mA$ be the matrix
of $\vA$ relative to $\vS$, and $\mA^\epsilon$ be the matrix of $\vA^\epsilon$ relative to $\vS^\epsilon$.
Let us summarize the properties of $\mA^\epsilon$:

\begin{proposition}
The matrix $\mA^\epsilon(z)=\mA(z)+\epsilon\mA'(z)$, where
\begin{equation}
\mA'=\begin{bmatrix}a'_{11}& a'_{12}\\a'_{21}& a'_{22}\end{bmatrix},\quad
\begin{array}{c}a'_{11},a'_{22}\in\Gamma(\p1,\vO(2))\\a'_{12}\in\Gamma(\p1,\vO(3))\\a'_{21}\in\Gamma(\p1,\vO(1)),\end{array}
\end{equation}
satisfies the following conditions:
\begin{enumerate}
\item For some $2\times 2$ matrix $\mS^\epsilon(z)=1+\epsilon\mS'(z)$, where the entries of $\mS'(z)$
are polynomials in $z$ (of arbitrary degree), we have $\mA^\epsilon(z)=\mS^\epsilon(z+1)^{-1}\mA(z)\mS^\epsilon(z)$.
\item For some $2\times 2$ matrix
\begin{equation*}
\mR^\epsilon(z)=\begin{bmatrix}1&0\\0&z^{-1}\end{bmatrix}(\mT(z^{-1})+\epsilon\mT'(z^{-1})),
\end{equation*}
where $\mT,\mT'$ are $2\times 2$ matrices over $\C[[z^{-1}]]$ and $\mT$ is invertible
(that is, $\det(\mT|_{z^{-1}=0})\ne0$), we have
\begin{equation*}
\mR^\epsilon(z+1)^{-1}\mA^\epsilon(z)\mR^\epsilon(z)=
\diag((\rho_1+\epsilon\rho'_1)(z^2+d_1z),(\rho_2+\epsilon\rho'_2)(z^2+d_2z)). \qed
\end{equation*}
\end{enumerate}
\label{pp:isomonodromic:existence}
\end{proposition}

Conversely, a matrix $\mA^\epsilon$ with such properties corresponds to the `continuous isomonodromy deformation'
of Theorem \ref{th:isomonodromic}(1). Actually, we can reformulate Theorem \ref{th:isomonodromic}(1)
(for $2n=4$, $\deg(\theta)=-1$) as
the following statement:

\begin{proposition}
Let $\mA(z)\in X_\theta$, $\theta\in\paramv_4$, $\deg(\theta)=-1$.
\begin{enumerate}
\item There is a deformation $\mA^\epsilon(z)$ that satisfies the conditions of
Proposition \ref{pp:isomonodromic:existence};

\item $\mA^\epsilon(z)$ is unique up to a d-gauge transformation
\begin{equation*}
\mA^\epsilon(z)\mapsto \mR^\epsilon(z+1)^{-1}\mA^\epsilon(z)\mR^\epsilon(z)
\end{equation*}
for a gauge matrix
$\mR^\epsilon(z)=1+\epsilon\mR'(z)$, where $\mR'(z)$ is of the form
\begin{equation*}
\mR'(z)=\begin{bmatrix}r'_{11}&r'_{12}\\0&
r'_{22}\end{bmatrix},\quad
\begin{array}{c}r'_{11},r'_{22}\in\C-\{0\}\\r'_{12}\in\Gamma(\p1,\vO(1)).\end{array}
\end{equation*}
\end{enumerate}
\qed
\end{proposition}

\subsection{Isomonodromy deformation of d-connections as $PVI$}
Let us now use coordinates $q,p$ on $M_\theta$ to write the connection of Theorem \ref{th:isomonodromic}(1)
as a system of differential equations on $p$ and $q$.
Suppose $(\vL,\vA)\in M_\theta$ and let $\mA\in X_\theta$ be the matrix of $\vA$ relative to
some trivialization $\vS:\vO\oplus\vO(-1)\iso\vL$. We need to find $\mA^\epsilon$ that satisfies the conditions
of Proposition \ref{pp:isomonodromic:existence}. As in Section \ref{sc:dPV},
it suffices to do so when $(\vL,\vA)$ belong to a dense subset of $M_\theta$; we can thus assume that
$q(\vL,\vA)\ne\infty$. We can then pick $\vS$ so that $\mA$ is of the form \eqref{eq:connmatV}.

We will look for $\mA^\epsilon$ in the form
\begin{equation*}
\mA^\epsilon(z)=\begin{bmatrix}a_{11}^\epsilon&a_{12}^\epsilon\\
a_{21}^\epsilon&a_{22}^\epsilon\end{bmatrix}=\mS^\epsilon(z+1)^{-1}\mA(z)\mS^\epsilon(z)
\end{equation*}
for the gauge matrix
\begin{equation*}
\mS^\epsilon(z)=1+\begin{bmatrix}s'_{11}&s'_{12}\\s'_{21}&s'_{22}\end{bmatrix}\epsilon\quad
\begin{array}{c}s'_{11},s'_{22}\in\Gamma(\p1,\vO(1))\\s'_{12}\in\Gamma(\p1,\vO(2))\\s'_{21}\in\C.\end{array}
\end{equation*}
(Actually, the proof of Proposition \ref{pp:isomonodromic:formal} shows that
$\mA^\epsilon$ is necessarily of this form.)
Then $\mA^\epsilon$ automatically satisfies Proposition \ref{pp:isomonodromic:existence}(1),
so we only need to make sure that Proposition \ref{pp:isomonodromic:existence}(2) is satisfied.
From Lemma \ref{lm:connmat} (which still holds
for d-connections that depend on $\epsilon$), we see that Proposition \ref{pp:isomonodromic:existence}(2) is
equivalent to the following equations:
\begin{gather*}
\det(\mA^\epsilon)=(z-a_1)(z-a_2)(z-a_3)(z-a_4)\rho^\epsilon_1\rho_2^\epsilon\\
a^\epsilon_{11}+a^\epsilon_{22}(1+z^{-1})=
(\rho^\epsilon_1+\rho^\epsilon_2)z^2+(d_1\rho^\epsilon_1+d_2\rho^\epsilon_2)z+t(z^{-1}),
\end{gather*}
where $\rho_i^\epsilon=\rho_i+\rho'_i\epsilon$, and $t(z^{-1})\in\C^\epsilon[[z^{-1}]]$ is a Taylor
series in $z^{-1}$ with coefficients in $\C^\epsilon$.
Solving these equations, we can find formulas for $q'$, $p'$ in terms of $\rho'_1$, $\rho'_2$, and $\theta$; here $q'$ and $p'$ are determined by the condition
\begin{equation*}
a^\epsilon_{21}(q+\epsilon q')=0\in\C^\epsilon,\qquad a^\epsilon_{11}(q+\epsilon q')=(p+\epsilon p')(q+\epsilon q'-a_3)(q+\epsilon q'-a_4).
\end{equation*}
The formulas for $q'$ and $p'$
can then be viewed as a system on differential equations on $q$ and $p$ (considered as functions of $\rho_i$):
\begin{equation}
\begin{split}
dq=&\dfrac{\rho_1 d\rho_2-\rho_2 d\rho_1}{\rho_1-\rho_2}
\left(\dfrac{p(q-a_3)(q-a_4)}{\rho_1\rho_2}-\dfrac{(q-a_1)(q-a_2)}{p}\right)
\cr
dp=&p\dfrac{d\rho_1-d\rho_2}{\rho_1-\rho_2}+
\dfrac{\rho_1 d\rho_2-\rho_2 d\rho_1}{\rho_1-\rho_2}\left(a_1+a_2-2q+\dfrac{p^2(a_3+a_4-2q)}{\rho_1\rho_2}\right.+
\cr
&\hphantom{p\dfrac{d\rho_1-d\rho_2}{\rho_1-\rho_2}+\dfrac{\rho_1 d\rho_2-\rho_2 d\rho_1}{\rho_1-\rho_2}}
\left.\dfrac{p}{\rho_1\rho_2}\left(d_1\rho_1+d_2\rho_2+
2q(\rho_1+\rho_2)\right)\right)
\end{split}
\label{eq:PVI}
\end{equation}

\begin{proof}[Proof of Theorem \ref{th:isomonodromic}(2)] We need to verify that \ref{eq:PVI} is obtained from the $PVI$
\ref{eq:standardPVI} by plugging in the formulas for $p^{PVI}$, $q^{PVI}$, $x_i$'s, and $\lambda_i^\pm$'s
(from Theorem \ref{th:dPV and PVI} and Section \ref{sc:Geometry of PVI}). This is a straightforward calculation.
\end{proof}

\begin{remark} Theorem \ref{th:isomonodromic}(2) can also be proved by an indirect argument. Indeed, both $PVI$
and \eqref{eq:PVI} define algebraic connections on the family $M\to P$ from Theorem \ref{th:isomonodromic}. The difference
between two such connections is a vector field on the moduli space $M_\theta$; on the other hand, it is known that
$M_\theta$ has no non-zero global vector fields (\cite[Theorem 3, Lemma 3]{AL0}, \cite[Proposition 2.1]{STT}).

Still another, more geometric, proof of Theorem \ref{th:isomonodromic}(2) uses the Mellin transform described
in Section \ref{sc:Fourier}. It is easy to see that under the transform, the continuous isomonodromy
deformation of d-connections (from Theorem \ref{th:isomonodromic}(1)) corresponds to the isomonodromy
deformation of ordinary connections, which is described by the sixth Painlev\'e equation.
\end{remark}

\subsection{Mellin transform}
\label{sc:Fourier}
In this section (which is completely independent from the rest of the paper), we sketch the geometric construction
underlying Theorem \ref{th:dPV and PVI}.
Fix $\theta\in\paramv_4$, $x\in X$, and $\lambda\in\Lambda$ as in Theorem \ref{th:dPV and PVI}.

Take $(\hat\vL,\nabla)\in M_{(x,\lambda)}$. For any $z\in\C$, consider the connection
\begin{equation*}
\nabla_z\eqd\nabla-z \zeta^{-1}d\zeta:\hat\vL\to\hat\vL\times\Omega_\p1(x_1+x_2+x_3+x_4),
\end{equation*}
where we denote by $\zeta$ the coordinate on $\p1$.
Recall that $x_1=0$, $x_4=\infty$, so subtraction of $z\zeta^{-1}d\zeta$ from $\nabla$ does not introduce new poles.
Denote by $\hat\vL_{*!}\supset\hat\vL$ the smallest quasi-coherent sheaf that contains $\hat\vL$ and such that
$\nabla_z(\hat\vL_{*!})\subset\hat\vL_{*!}$ for all $z\in\C$.
(In terms of $D$-modules, $\hat\vL_{*!}$ can be constructed by taking the intermediate extension of $(\hat\vL,\nabla_z)$ from $\p1-\{x_1,x_2,x_3,x_4\}$
to $\p1-\{0,\infty\}$ and then extending to $\p1$.) Consider the first de Rham cohomology group $H^1_{DR}(\hat\vL_{*!},\nabla_z)$.
Since $\hat\vL_{*!}$ and $\hat\vL_{*!}\otimes\Omega_\p1$ have no higher cohomologies, it can be computed by the formula
\begin{equation*}
H^1_{DR}(\hat\vL_{*!},\nabla_z)=\coker(\nabla_z:\Gamma(\p1,\hat\vL_{*!})\to\Gamma(\p1,\hat\vL_{*!}\otimes\Omega_\p1)).
\end{equation*}
$H^1_{DR}(\hat\vL_{*!},\nabla_z)$ depends on $z$ in an algebraic way; more precisely, it is the fiber over $z\in\C$
of a natural quasicoherent sheaf $\vL_{*!}$ on $\p1-\{\infty\}$. The sheaf $\vL_{*!}$ is the Mellin transform of
$\hat\vL_{*!}$ in terms of \cite{L}.

Consider now the rational map $a:\hat\vL\isorat\hat\vL:s\mapsto \zeta s$. Note that $a$ satisfies the relation
$a\circ\nabla_z=\nabla_{z+1}\circ a$. It is also easy to see that $a$ induces an automorphism of $\hat\vL_{*!}$;
therefore, it becomes an isomorphism
of $D$-modules (that is, quasicoherent sheaves with connections) $(\hat\vL_{*!},\nabla_z)\iso(\hat\vL_{*!},\nabla_{z+1})$. Hence $a$ yields an identification.
\begin{equation*}
\widetilde\vA(z):H^1_{DR}(\hat\vL_{*!},\nabla_z)\iso H^1_{DR}(\hat\vL_{*!},\nabla_{z+1}).
\end{equation*}
As $z\in\C$ varies, we can view $\widetilde\vA(z)$ as a d-connection on the quasicoherent sheaf $\vL_{*!}$. One can check that $\vL_{*!}$ contains a unique
coherent locally free subsheaf of rank $2$ (that is, a rank $2$ vector bundle) $\vL\subset\vL_{*!}$ such that
\begin{equation*}
\vA(z)\eqd (z-a_3)(z-a_4)\widetilde\vA(z)
\end{equation*}
is a d-connection
of type $\theta$ on $\vL$. The correspondence
\begin{equation*}
(\hat\vL,\nabla)\mapsto(\vL,\vA)
\end{equation*}
gives a map $M_{(x,\lambda)}\to M_\theta$. Note that the scalar multiple $(z-a_3)(z-a_4)$ also appears in Remark \ref{rm:dPV and Gamma}.

To describe the inverse map $M_\theta\to M_{(x,\lambda)}$, let us reconstruct $(\hat\vL,\nabla)$ from $(\vL,\vA)$. For any $\zeta\in\C-\{0\}$,
consider the d-connection
\begin{equation*}
\widetilde\vA_\zeta\eqd\zeta^{-1}\dfrac{\vA}{(z-a_3)(z-a_4)}
\end{equation*}
on $\vL$. Let $\vL_{*!}$ be the smallest quasicoherent sheaf on $\p1$ that contains $\vL$ and such that $\widetilde\vA_\zeta$ induces an isomorphism
$(\vL_{*!})_z\to(\vL_{*!})_{z+1}$ for all $z$ and $\zeta$ (the quotient $\vL_{*!}/\vL$ is the direct sum of length $1$ skyscraper sheaves
supported at points $a_1,a_1-1,a_1-2,\dots;a_2,a_2-1,\dots;a_3+1,a_3+2,\dots,a_4+1,a_4+2,\dots$). For any $\zeta\in\C-\{0\}$,
we obtain a structure of a $\Z$-equivariant sheaf on $\vL_{*!}$,
where $1\in\Z$ acts on $\p1$ by $z\mapsto z+1$ and on $\vL_{*!}$ by $\widetilde\vA_\zeta$ (in some sense, $\vL_{*!}$ is obtained from $\vL$ by
an `intermediate extension' for $\Z$-equivariant sheaves). Consider the corresponding equivariant cohomology
group $H^1_\Z(\vL_{*!},\widetilde\vA_\zeta)$, which can be computed by the formula
\begin{equation*}
H^1_\Z(\vL_{*!},\widetilde\vA_\zeta)=\coker(\widetilde\vA_\zeta-1:\Gamma(\p1,\vL_{*!})\to\Gamma(\p1,\vL_{*!})).
\end{equation*}
$H^1_\Z(\vL_{*!},\widetilde\vA_\zeta)$ is the fiber over $\zeta\in\C-\{0\}$
of the quasicoherent sheaf $\hat\vL_{*!}$ on $\p1-\{\infty,0\}$.

For every $\zeta\in\C-\{0\}$, consider the rational map
\begin{equation*}
\delta(\zeta):\vL\dashrightarrow\vL:s\mapsto z\zeta^{-1}s.
\end{equation*}
$\delta(\zeta)$ induces a regular map $\vL_{*!}\to\vL_{*!}$, and, therefore, a map
\begin{equation*}
\delta_*(\zeta):\Gamma(\p1,\vL_{*!})\to\Gamma(\p1,\vL_{*!}).
\end{equation*}
The map $\delta_*(\zeta)$ satisfies the following commutativity relation
\begin{equation*}
\delta_*(\zeta)\widetilde A_\zeta=\widetilde A_\zeta\delta_*(\zeta)+\dfrac{d\widetilde A_\zeta}{d\zeta}.
\end{equation*}
Now let us consider the trivial quasicoherent sheaf over $\p1-\{0,\infty\}$ whose fiber over every point $\zeta\in\p1-\{0,\infty\}$ equals $\Gamma(\p1,\vL_{*!})$.
The formula $\widetilde\vA_\zeta-1$ gives an endomorphism of this sheaf; the cokernel of the endomorphism is $\hat\vL_{*!}$. Notice now that
$\widetilde\vA_\zeta-1$ is horizontal with respect to the connection $\nabla=d+\delta_*(\zeta)d\zeta$ on the sheaf. Therefore, $\nabla$ induces a connection
$\hat\vL_{*!}\to\hat\vL_{*!}\otimes\Omega_\p1$ (which we will also denote by $\nabla$).
Finally, $\hat\vL\subset\hat\vL_{*!}$ can be reconstructed as the only coherent locally free subsheaf of rank $2$ such that
$\nabla$ is a connection of type $(x,\lambda)$ on $\hat\vL$.

\section{Difference $PVI$}
In this section, we study $M_\theta$ for $\theta\in\paramvi_k$. We will need suitable versions of several statements from
Section \ref{sc:General d-connections}.

\subsection{}
\begin{proposition}[cf. Proposition \ref{pp:infty}] Suppose that the matrix
$\mA(z)=\sum\limits_{i\le n}\mA_i z^i$ over $\C((z^{-1}))$
satisfies the following condition:
\begin{equation}
\label{cn:dpvi} \begin{array}{l} \text{The leading term $\mA_n$ is a non-zero scalar matrix}\cr
\text{while all eigenvalues of the next term $\mA_{n-1}$ are distinct.}
\end{array}
\end{equation}
Then there exists a gauge matrix
$\mR(z)=\sum\limits_{i\le0}\mR_iz^i$ with invertible $\mR_0$ such
that
\begin{equation}
\mR(z+1)^{-1}\mA(z)\mR(z)=\mA'_n z^n+\mA'_{n-1} z^{n-1},
\end{equation}
where $\mA'_n$ and $\mA'_{n-1}$ is diagonal. $\mR(z)$ is uniquely determined up
to right multiplication by a
permutation matrix and a constant diagonal matrix. \qed
\label{pp:inftyvi}
\end{proposition}

As before, we will denote the only eigenvalue of $\mA'_n$ by $\rho=\rho_1=\dots=\rho_n$, and the eigenvalues
of $\mA'_{n-1}$ by $\rho d_1,\dots,\rho d_n$. It is easy to see that $\mA_n=\mA'_n$ (so $\rho$ is also the eigenvalue
of $\mA_n$) and $\mA_{n-1}$ is conjugate to $\mA'_{n-1}$ (so $\rho d_1,\dots,\rho d_n$ are also eigenvalues of $\mA'_{n-1}$;
this can be thought of as a version of Remark \ref{rm:characteristic polynomial}.

\begin{proposition}[cf. Proposition \ref{pp:modv}] Suppose $\theta=(a_1,\dots,a_k;\rho,\rho,d_1,d_2;n)$,
and $d_1\ne d_2$. Let $(\vL',\vA')$ be an elementary upper
modification of $(\vL,\vA)\in M_\theta$ given by $(x\in\p1;l\subset\vL_x)$. Then the
only cases when $(\vL',\vA')$ belongs to $M_{\theta'}$ for some
$\theta'\in\param$ are as follows:
\begin{enumerate}
\item If $x=\infty$, then $l$ must be an
eigenspace of $\vA_{n-1}:\vL_\infty\to\vL_\infty$ (the second term of
$\vA=\rho z^n+\vA_{n-1} z^{n-1}+\text{lower order terms}$). If, for instance,
$l=\ker(\vA_{n-1}-\rho d_1)\subset\vL_{\infty}$, then
$\theta'=(a_1,\dots,a_k;\rho_1,\rho_2,d_1-1,d_2;n)$,
and an analogous formula holds when $l=\ker(\vA_{n-1}-\rho d_2)$.

\item If $x=a_i$ is a zero of $\vA$ and $x-1\ne
a_j$ is not, then $l$ must be the kernel of
$\vA(x):\vL_x\to\vL_{x+1}$; in this case,
$\theta'=(a_1,\dots,a_i-1,\dots,a_k;\rho_1,\rho_2,d_1,d_2;n)$.
\end{enumerate}

In either case, the elementary modifications define an isomorphisms
$M_\theta\iso M_{\theta'}$.
\qed
\label{pp:modvi}
\end{proposition}

\begin{corollary}
Suppose $\theta\in\Theta_k$ satisfies \eqref{cn:mod}, \eqref{cn:rhoalt}. Then
$M_\theta$ is naturally isomorphic to $M_{\theta'}$ whenever $\theta'$ is
obtained from $\theta$ by adding integers to $a_i$'s and $d_i$'s. \qed
\end{corollary}

\begin{lemma}[cf. Corollary \ref{co:height}]
Suppose $(\vL,\vA)\in M_\theta$ and suppose that $\theta\in\param_{2n}$ satisfies \eqref{cn:irred}, \eqref{cn:rhoalt}.
If $\vL\simeq\vO(n_1)\oplus\vO(n_2)$, then $|n_1-n_2|\le n-1$.
\label{lm:heightvi}
\end{lemma}
\begin{proof} The proof repeats that of Corollary \ref{co:height}; the only difference is that the order of pole of
$\alpha$ at $\infty$ cannot exceed $n-1$ (because the coefficient of $z^n$ in $\alpha$
is an off-diagonal element of a scalar matrix; that is, zero).
\end{proof}

\subsection{Proof of Theorems \ref{th:Geometry of dPVI}, \ref{th:dPVI}}
The proof of Theorem \ref{th:Geometry of dPVI} follows the
same ideas as that of Theorem \ref{th:Geometry of dPV}.
Fix $\theta\in\paramvi_6$, $\deg(\theta)=-1$.
For any $(\vL,\vA)\in M_\theta$, Lemma \ref{lm:heightvi} implies $\vL\simeq\vO\oplus\vO(-1)$. Choosing an isomorphism
$\vS:\vO\oplus\vO(-1)\iso\vL$, we can write $\vA$ as a matrix
\begin{equation}
\mA=\begin{bmatrix}a_{11}& a_{12}\\a_{21}& a_{22}\end{bmatrix},\quad
\begin{array}{c}a_{11},a_{22}\in\Gamma(\p1,\vO(3))\\a_{12}\in\Gamma(\p1,\vO(4))\\a_{21}\in\Gamma(\p1,\vO(2)).\end{array}
\label{eq:connmatvi}
\end{equation}
Choosing a different isomorphism $\vS$ replaces $\mA$ with its d-gauge transformation \eqref{eq:gauge2}, where the gauge matrix
$\mR$ is given by \eqref{eq:gaugemat}.

\begin{lemma}[cf. Lemma \ref{lm:connmat}]
Let $\vA$ be a d-connection on $\vO\oplus\vO(-1)$; its matrix $\mA$ is of the form \eqref{eq:connmatvi}.
We claim that $\vA$ is of type $\theta$ if and only if $\mA$ satisfies the following conditions:
\begin{gather*}
a_{12}\in\Gamma(\p1,\vO(3));\qquad a_{21}\in\Gamma(\p1,\vO(1));\qquad
a_{11}-\rho z^3,a_{22}-\rho z^3\in\Gamma(\p1,\vO(2))\\
\det(\mA)=(z-a_1)(z-a_2)(z-a_3)(z-a_4)(z-a_5)(z-a_6)\rho^2\\
(a_{11}-\rho z^3)(a_{22}(1+z^{-1})-\rho z^3)-a_{12}a_{21}=d_1d_2\rho^2 z^4+\text{lower order terms}.
\end{gather*}
\label{lm:connmatvi}\qed
\end{lemma}

\begin{remark*}
The last condition of the lemma can be more naturally written as
\begin{equation*}
\det(\mR(z+1)^{-1}\mA\mR(z)-\rho z^3)=d_1d_2\rho^2 z^4+\text{lower order terms},
\end{equation*}
where $\mR(z)\eqd\diag(1,z^{-1})$ is a trivialization of $\vO\oplus\vO(-1)$ near $\infty\in\p1$.
\end{remark*}

We can now think of $M_\theta$ as the quotient
of the space of all matrices \eqref{eq:connmatvi} that satisfy Lemma \ref{lm:connmatvi} modulo d-gauge transformations
with gauge matrices \eqref{eq:gaugemat} (cf. Corollary \ref{co:dPV as quotient}).
For any matrix \eqref{eq:connmatvi} that satisfies Lemma \ref{lm:connmatvi}, denote by $q\in\p1$ the only zero of $a_{21}$, and
set $\tilde p\in(\vO(3))_q$. It is easy to see that $q$ and $\tilde p$ do not change
under d-gauge transformations with gauge matrices \eqref{eq:gaugemat}; therefore, $\tilde P\eqd(q,\tilde p)$ can be viewed as a map
$M_\theta\to\tK$, where $\tK\eqd\V(\vO(3)^\vee)$ is the total space of the line bundle $\vO(3)$. We can now use the map $\tilde P$
for a geometric description of $M_\theta$ (we are using the notation of Theorem \ref{th:Geometry of dPV tilde}):

\begin{theorem}
\begin{enumerate}
\item The map $\tilde P:M_\theta\to\tK$ is a regular birational
morphism of smooth algebraic surfaces.

\item Let $\tilde\sigma_1:\tK_1\to\tK$ be the blow-up of $\tK$ at
the following $7$ points: $(a_i,0(a_i))$ ($i=1,\dots,6$) and
$(\infty,(\rho z^3)(\infty))$. Let $\sigma_2:\tK_2\to
\tK_1$ be the blow-up of $\tK_1$ at the two points
$(\infty,(\rho z^3+\rho d_jz^2)'(\infty))$, $j=1,2$ (these points
belong to the preimage $\tilde\sigma_1^{-1}(\infty,(\rho z^3)(\infty))\subset\tK_1$).
Then the map $\tP$ induces an open embedding
$\tP_2:M_\theta\hookrightarrow \tK_2$.

\item The complement to $\tP_2(M_\theta)$ in $\tK_2$ is the union
of the proper preimages of the following curves: the zero section
$\{(z,0(z)):z\in\p1\}\subset\tK$, the fiber at infinity
$\{(\infty,az^3(\infty)):a\in\C\}\subset\tK$, and the
exceptional curve
$\tilde\sigma_1^{-1}(\infty,(\rho z^3)(\infty))\subset\tK_1$.
\end{enumerate}
\label{th:Geometry of dPVI tilde}\qed
\end{theorem}

The proof of Theorem \ref{th:Geometry of dPVI tilde} is completely analogous to that of Theorem \ref{th:Geometry of dPV tilde}
(Section \ref{sc:Geometry of dPV tilde}). Now Theorem \ref{th:Geometry of dPVI} easily follows: we set
\begin{equation*}
p\eqd\dfrac{\tilde p}{(q-a_4)(q-a_5)(q-a_6)},
\end{equation*}
and it is not hard to check that the map $P\eqd(q,p):M_\theta\to(\p1)^2$ (which is birational by Theorem \ref{th:Geometry of dPVI tilde})
is regular and induces an embedding $M_\theta\hookrightarrow K_2$ with the required properties.

\begin{proof}[Proof of Theorem \ref{th:dPVI}]
The proof repeats the proof of Theorem \ref{th:dPV} (given in Section \ref{sc:dPV}) almost word-for-word (of course, the calculations involved are somewhat more complicated).
The only real difference is formulas \eqref{eq:connmatV}, \eqref{eq:connmatV'}; the corresponding formulas in our case are
\begin{align*}
\mA&=\begin{bmatrix}z^3-q^3+p(q-a_4)(q-a_5)(q-a_6)& a_{12}\\z-q&a_{22}\end{bmatrix},\quad
a_{22},a_{12}\in\Gamma(\p1,\vO(3));\\
\mA'&=\begin{bmatrix}z^3-(q')^3+p'(q'-a_3)(q'-a_4)(q'-a_6)& a'_{12}\\z-q'&a'_{22}\end{bmatrix},\quad
a'_{22},a'_{12}\in\Gamma(\p1,\vO(3)).
\end{align*}
\end{proof}

\subsection{Degeneration to difference $PV$}
\label{sc:dPVI and dPV}
Given
\begin{equation*}
\tilde\theta=(\tilde a_1,\tilde a_2,\tilde a_3,\tilde a_4;\tilde \rho_1,\tilde \rho_2,\tilde d_1,\tilde d_2;2)
\in\paramv_4,
\end{equation*}
let us define $\theta(t)$ for $t\in\C-\{0\}$ by
\begin{equation*}
\theta(t)=(\tilde a_1,\tilde a_2,-\tilde \rho_1/t,-\tilde \rho_2/t,\tilde a_3,\tilde a_4;1,1,\tilde d_1+
(\tilde\rho_1/t),\tilde d_2+(\tilde \rho_2/t);3);
\end{equation*}
clearly, $\theta(t)\in\paramvi_6$ for all but countably many $t$. Denote the components of
$\theta=\theta(t)$ by $a_i=a_i(t)$, $d_j=d_j(t)$. Formulas \eqref{eq:dPVI} define a family
of equations depending on parameter $t\in\C-\{0\}$. Let us show that the difference $PV$ \eqref{eq:dPV}
is the limit of this family as $t\to 0$.

Replace $p$ with a new variable $\tilde p\eqd (\tilde\rho_2+qt)p$; accordingly, set
$\tilde p'\eqd (\tilde\rho_2+q't)p'$. After we plug the formulas for $\theta(t)$, $\tilde p$, and $\tilde p'$ into
\eqref{eq:dPVI}, it becomes the following system:
\begin{equation}
\begin{cases}
q+q'=\tilde a_3+\tilde a_4+\dfrac{\tilde\rho_1(\tilde d_1+\tilde a_4+\tilde a_5)}{\tilde p-\tilde\rho_1}
+\dfrac{\tilde\rho_2(\tilde d_2+\tilde a_4+\tilde a_5)}{\tilde p-\tilde\rho_2}+O(t)
\cr
\tilde p\tilde p'=\dfrac{(q'-\tilde a_1+1)(q'-\tilde a_2+1)}{(q'-\tilde a_3)(q'-\tilde a_4)}\cdot\tilde\rho_1\tilde\rho_2+
O(t),
\end{cases}
\end{equation}
where $O(t)$ stands for a Taylor series in $t$ with no constant term.
This is exactly the difference $PV$ equation \eqref{eq:dPV}.

\begin{remark}
\label{re:degenerate geometry}
The degeneration of \eqref{eq:dPVI} to \eqref{eq:dPV} has a clear geometric meaning; let us sketch it.
It is easy to
construct a family of moduli spaces $\nu:N\to\A1$ such that the fiber $\nu^{-1}(t)$ over $t\in\A1-\{0\}$ equals
$M_{\theta(t)}$ whenever $\theta(t)\in\paramvi_6$, while $\nu^{-1}(0)=M_{\tilde\theta}$. Similarly,
one can define a family $\nu':N'\to\A1$ such that $(\nu')^{-1}(t)=M_{\theta'(t)}$ if
$t\ne 0$, $\theta(t)\in\paramvi_6$ and that $(\nu')^{-1}(0)=M_{\tilde\theta'}$. Here
\begin{equation*}
\tilde\theta'=(\tilde a_1+1,\tilde a_2+1,\tilde a_3,\tilde a_4;\tilde\rho_1,\tilde\rho_2,\tilde d_1-1,\tilde d_2-1;2),
\end{equation*}
\begin{equation*}
\theta'(t)=(\tilde a_1+1,\tilde a_2+1,-\tilde \rho_1/t,-\tilde \rho_2/t,\tilde a_3,\tilde a_4;1,1,\tilde d_1+
(\tilde\rho_1/t)-1,\tilde d_2+(\tilde \rho_2/t)-1;3).
\end{equation*}
The modification of d-connections defines a rational isomorphism $N\isorat N'$ that is regular over a neighborhood
of $0\in\A1$; this isomorphism is given by \eqref{eq:dPVI} if $t\ne0$ and $\theta(t)\in\paramvi_6$ and by
\eqref{eq:dPV} if $t=0$.
\end{remark}

\subsection{Degeneration to classical $PVI$}
\label{sc:dPVI and PVI}
Let us now show how difference $PVI$ \eqref{eq:dPVI} degenerates into the classical $PVI$. Fix
\begin{equation*}
\tilde\theta=(\tilde a_1,\tilde a_2,\tilde a_3,\tilde a_4;\tilde \rho_1,\tilde \rho_2,\tilde d_1,\tilde d_2;2)
\in\paramv_4,
\end{equation*}
and set
\begin{equation*}
\theta(t)\eqd(-\tilde \rho_1/t,-\tilde \rho_2/t,\tilde a_1,\tilde a_2,\tilde a_3,\tilde a_4;1,1,\tilde d_1+
(\tilde\rho_1/t),\tilde d_2+(\tilde \rho_2/t);3)\quad(t\in\C-\{0\});
\end{equation*}
again, $\theta(t)\in\paramvi_6$ for all but countably many $t$. Let us also set
\begin{equation*}
\theta'(t)\eqd(-(\tilde \rho_1/t)-1,-(\tilde \rho_2/t)-1,\tilde a_1,\tilde a_2,\tilde a_3,\tilde a_4;1,1,\tilde d_1+
(\tilde\rho_1/t)+1,\tilde d_2+(\tilde \rho_2/t)+1;3),
\end{equation*}
so that $dPVI$ is an isomorphism $M_{\theta(t)}\iso M_{\theta'(t)}$. Note that the formula for $\theta'(t)$
is obtained from the formula for $\theta(t)$ if we substitute
\begin{equation}
\tilde\rho'_i\eqd\tilde\rho_i+t,
\label{eq:rho'}
\end{equation}
for $\tilde\rho_i$, $i=1,2$.

Let us replace $p$ with $\tilde p\eqd (q-\tilde a_2)tp$; accordingly, set $\tilde p'\eqd (q'-\tilde a_2)tp'$.
Then \eqref{eq:dPVI} can be written as
\begin{equation}
\begin{cases}
\dfrac{q'-q}{t}=\dfrac{(q-\tilde a_3)(q-\tilde a_4)}{\tilde\rho_1\tilde\rho_2}\tilde p-
\dfrac{(q-\tilde a_1)(q-\tilde a_2)}{\tilde p}+O(t)\cr
\dfrac{\tilde p'-\tilde p}{t}=\tilde a_1+\tilde a_2-2q+\dfrac{2(\tilde\rho_1+\tilde\rho_2)q+
\tilde d_1\tilde\rho_1+\tilde d_2\tilde\rho_2}{\tilde\rho_1\tilde\rho_2}\tilde p+
\dfrac{\tilde a_3+\tilde a_4-2q}{\tilde\rho_1\tilde\rho_2}\tilde p^2+O(t),
\end{cases}
\label{eq:dPVIlim}
\end{equation}
where $(q,p)$ are the coordinates on $M_{\theta(t)}$ and $(q',p')$ are the coordinates on $M_{\theta'(t)}$.
As $t\to 0$, the left hand sides tend to derivatives of $q$ and $p$ with respect to $t$. Similarly,
\eqref{eq:rho'} becomes the expression
\begin{equation*}
\dfrac{d\tilde\rho_i}{d t}=1, \quad (i=1,2);
\label{eq:rhoder}
\end{equation*}
all other parameters $\tilde a_1,\dots,\tilde a_4;\tilde d_1,\tilde d_2$ do not depend on $t$.
Now it is easy to see that \eqref{eq:dPVIlim} is obtained from \eqref{eq:PVI}
(which is equivalent to the sixth Painlev\'e equation) by changing variables from $\tilde\rho_1$, $\tilde\rho_2$ to $t$.

The degeneration of \eqref{eq:dPVI} to \eqref{eq:dPVIlim} has a geometric interpretation similar to that given
for the degeneration to \eqref{eq:dPV} (Remark \ref{re:degenerate geometry}). The details are left to the reader.
\nocite{*}
\bibliographystyle{plain}
\bibliography{dPVIrefs}
\end{document}